\documentclass[11pt]{amsart}

\usepackage[english]{babel}
\usepackage[T1]{fontenc}
\usepackage[utf8]{inputenc}

\usepackage{mathpazo}
\usepackage{graphicx}
\usepackage{microtype}
\usepackage{enumerate}
\usepackage{csquotes}
\usepackage[dvipsnames]{xcolor}
\usepackage[bookmarks, colorlinks, citecolor=MidnightBlue!75!black, urlcolor=MidnightBlue!75!black, linkcolor=MidnightBlue!75!black]{hyperref}

\calclayout

\AtBeginDocument{%
   \def\MR#1{}
}

\usepackage{amsmath, amssymb, amsthm, amscd}
\usepackage{commath, mathtools, mathrsfs, thmtools}
\usepackage{tikz, tikz-cd}
\usetikzlibrary{matrix, arrows, decorations.pathmorphing}

\theoremstyle{plain}
\newtheorem{theorem}{Theorem}[section]
\newtheorem{lemma}[theorem]{Lemma}
\newtheorem{proposition}[theorem]{Proposition}
\newtheorem{corollary}[theorem]{Corollary}

\newtheorem*{TA}{Theorem A}
\newtheorem*{TB}{Theorem B}
\newtheorem*{TC}{Theorem C}
\newtheorem*{TD}{Theorem D}

\theoremstyle{definition}
\newtheorem{remark}[theorem]{Remark}
\newtheorem{definition}[theorem]{Definition}
\newtheorem{example}[theorem]{Example}

\newtheorem*{conjecture*}{Conjecture}
\newtheorem*{remark*}{Remark}
\newtheorem*{definition*}{Definition}

\newcommand{\on}{\operatorname}

\newcommand{\aff}{\on{Aff}}
\newcommand{\affgrp}{\on{AffGrp}_{k}}
\newcommand{\aut}{\on{Aut}}

\newcommand{\C}{\mathbb{C}}

\newcommand{\fgrp}{\on{FGrp}}
\newcommand{\FLi}{\on{FL}_{\infty}}

\newcommand{\gal}{\on{Gal}}

\newcommand{\hz}{\hat{\mathbb{Z}}}
\newcommand{\id}{\on{id}}

\newcommand{\isom}{\on{Isom}}

\newcommand{\mc}{\mathcal}
\newcommand{\mf}{\mathfrak}

\newcommand{\op}{\on{op}}

\newcommand{\Pic}{\on{Pic}}
\newcommand{\Q}{\mathbb{Q}}

\DeclareRobustCommand{\rchi}{{\mathpalette\irchi\relax}}
\newcommand{\irchi}[2]{\raisebox{0.9\depth}{$#1\chi$}}
\newcommand{\rra}{\rightrightarrows}

\newcommand{\s}{\subseteq}

\newcommand{\spec}{\on{Spec}}

\newcommand{\uaut}{\underline{\on{Aut}}}
\newcommand{\uisom}{\underline{\on{Isom}}}
\newcommand{\upi}{\underline{\pi}}

\newcommand{\xar}{\xrightarrow}
\newcommand{\Z}{\mathbb{Z}}

\renewcommand{\epsilon}{\varepsilon}
\renewcommand{\H}{\on{H}}
\renewcommand{\hom}{\on{Hom}}
\renewcommand{\injlim}{\varinjlim}

\renewcommand{\O}{\mc{O}}

\renewcommand{\phi}{\varphi}
\renewcommand{\projlim}{\varprojlim}
\renewcommand{\set}{\on{Set}}
\renewcommand{\tilde}{\widetilde}
\renewcommand{\hat}{\widehat}

\author{Giulio Bresciani}
\address{Freie Universit\"at Berlin, Arnimallee 3, 14195, Berlin, Germany}
\email{gbresciani@zedat.fu-berlin.de}
\title{Some implications between Grothendieck's anabelian conjectures}
\date{}

\thanks{The author is supported by the DFG Priority Program "Homotopy Theory and Algebraic Geometry" SPP 1786}

\begin{document}

\begin{abstract}
Grothendieck gave two forms of his "main conjecture of anabelian geometry", i.e. the section conjecture and the hom conjecture. He stated that these two forms are equivalent and that if they hold for hyperbolic curves then they hold for elementary anabelian varieties too. We state a stronger form of Grothendieck's conjecture (equivalent in the case of curves) and prove that Grothendieck's statements hold for our form of the conjecture. We work with DM stacks, rather than schemes. If $X$ is a DM stack over $k\s\C$ finitely generated over $\Q$, we prove that whether $X$ satisfies the conjecture or not depends only on $X_{\C}$. We prove that the section conjecture for hyperbolic orbicurves stated by Borne and Emsalem follows from the conjecture for hyperbolic curves.
\end{abstract}

\maketitle
\tableofcontents

\section{Introduction}

\subsection{The main conjecture of anabelian geometry}

In his letter to Faltings \cite{gro97}, Grothendieck gave two forms of his "main conjecture of anabelian geometry", the \emph{hom conjecture} and the \emph{section conjecture}, for \emph{anabelian} varieties. He refrained from defining precisely the class of anabelian varieties, but he said that it certainly contained smooth, hyperbolic curves and the so-called elementary anabelian varieties, i.e. those obtained by subsequent fibrations from smooth hyperbolic curves. He said that being anabelian is a purely geometric property, i.e. whether $X/k$ is anabelian depends only on $X_{\bar{k}}$, or $X_{\C}$ is $k\s\C$, that the section form implies the hom form and that, if the main conjecture holds for proper, hyperbolic curves, then it holds for proper, elementary anabelian varieties. Up to our knowledge, there is no proof of these statements in the literature. 

Let us recall the two forms of the main conjecture. If $G,H$ are extensions of a group $\Gamma$, then $\on{Hom-ext}_{\Gamma}(G,H)$ is the set of homomorphisms $G\to H$ which commute with the projection to $\Gamma$ modulo the natural action of $\ker(H\to\Gamma)$ by conjugation. If $k$ is a field, we denote by $\Gamma_{k}=\gal(k^{s}/k)$ the absolute Galois group. If $T,X$ are geometrically connected over $k$ with étale fundamental groups $\pi_{1}(T),\pi_{1}(X)$ (we omit base points), there is a natural map
\[\hom_{k}(T,X)\to\on{Hom-ext}_{\Gamma_{k}}(\pi_{1}(T),\pi_{1}(X)).\]

\begin{conjecture*}[Grothendieck, "Hom conjecture"]\label{gc1}
	Let $k$ be finitely generated over $\Q$. If $T/k$ is a smooth variety and $X/k$ is a smooth, proper anabelian variety, then
	\[\hom_{k}(T,X)\to\on{Hom-ext}_{\Gamma_{k}}(\pi_{1}(T),\pi_{1}(X))\]
	is a bijection.
\end{conjecture*}

There is a weaker form of the hom conjecture which restricts the attention to dominant morphisms, and this weaker form has been famously proved by Mochizuki for hyperbolic curves in \cite{moc99}. We stress that the hom conjecture is strictly stronger: Mochizuki's result only applies to \emph{open} homomorphisms $\pi_{1}(T)\to\pi_{1}(X)$, while the hom conjecture regards all homomorphisms.

The second form of the main conjecture is the so called \emph{section conjecture}, which is just the hom conjecture for $T=\spec k$.

\begin{conjecture*}[Grothendieck, "Section conjecture"]\label{gc2}
	Let $k$ be finitely generated over $\Q$. If $X$ is a smooth, proper anabelian variety over $k$, then
	\[X(k)\to\on{Hom-ext}_{\Gamma_{k}}(\Gamma_{k},\pi_{1}(X))\]
	is a bijection.
\end{conjecture*}

\subsection{A strengthening of the main conjecture}

The set $\on{Hom-ext}_{\Gamma_{k}}(\Gamma_{k},\pi_{1}(X))$ has a natural interpretation as the set of isomorphism classes of the category of rational points of the \emph{étale fundamental gerbe} $\Pi_{X/k}$, see \cite[§ 9]{bv15}. If $X$ is a DM stack, $X(k)$ has a natural structure of category too (rather than just a set) and Grothendieck's section map extends naturally to a functor
\[X(k)\to\Pi_{X/k}(k).\]
It is then natural to ask for an equivalence of categories rather than a mere bijection. In fact, already Grothendieck had pointed out that for stacks (which he called "multiplicities") the correct statement needs an equivalence of categories, see \cite[pg. 7]{gro97}. For hyperbolic curves the category structure is known to be trivial (i.e. it is just a set) on both sides.

\begin{definition*}
	Let $X$ be a smooth, proper, geometrically connected Deligne-Mum\-ford stack over a field $k$ of characteristic $0$. We say that $X$ is \emph{printable} (resp. \emph{fundamentally fully faithful, or fff}) if the natural morphism
	\[X(k')\to\Pi_{X/k}(k')\]
	is an equivalence (resp. fully faithful) for every finitely generated extension $k'/k$. 
\end{definition*}

The name "printable" is meant to suggest that $\Pi_{X/k}$ "prints" $X$, similarly to how an algebraic space represents a sheaf.

For smooth, proper, hyperbolic curves over a field $k$ finitely generated over $\Q$, printability is equivalent to the section conjecture over all finitely generated extensions of $k$, see \autoref{curvessect}.

\subsection{Results of the paper}

We prove that printability is a geometric property.

\begin{TA}[\ref{anabext}]
	Let $k'/k$ be a finitely generated extension, and $X$ a smooth, proper, geometrically connected DM stack over $k$. Then $X$ is printable (resp. fff) if and only if $X_{k'}$ is printable (resp. fff).
	
	As a consequence, if $k\s\C$ is finitely generated over $\Q$, whether or not $X$ is printable depends only on $X_{\C}$.
\end{TA}

We remark that in our proof of Theorem A the categorical structure is crucial, even for schemes: we do not know whether the same result holds if we replace printability with the analogous statement asking only for a bijection.

The first non-trivial example of expected anabelian DM stacks are hyperbolic orbicurves, see \cite{be14}. There is a natural notion of rational Euler characteristic for orbicurves, and hyperbolic ones are those with negative characteristic.

\begin{conjecture*}[Borne, Emsalem]
	Smooth, proper, hyperbolic orbicurves over finitely generated extensions of $\Q$ satisfy the section conjecture.
\end{conjecture*}

We prove that hyperbolic orbicurves are fff (in particular, they satisfy the injectivity part of the conjecture), and that the section conjecture for them is equivalent to the one for hyperbolic curves.

\begin{TB}[\ref{eq}]
	Let $k$ be finitely generated over $\Q$.
	\begin{itemize}
		\item A smooth, proper orbicurve is fff if and only if its Euler characteristic is less than or equal to $0$.
		\item If smooth, proper, hyperbolic curves satisfy the section conjecture over every finite extension $k'/k$, the section map is an equivalence for smooth, proper, hyperbolic orbicurves over $k$.
		\item Smooth, proper, hyperbolic orbicurves are printable if and only smooth, proper, hyperbolic curves are printable.
	\end{itemize}
\end{TB}

Thanks to Theorem B and an idea of N. Borne and M. Emsalem, we give a new, natural proof of the fact that the section conjecture for proper curves implies the section conjecture for affine curves using orbicurves as an intermediate step, see \autoref{propen}.

We then show that the section conjecture implies the hom conjecture.

\begin{TC}[\ref{point}]
	Let $X$ be a smooth, proper, geometrically connected DM stack and $T$ a locally noetherian, normal scheme over $k$. Assume that, for every $t\in T$, the residue field $k(t)$ is finitely generated over $k$. If $X$ is fff, then $X(T)\to\Pi_{X}(T)$ is fully faithful. If $X$ is printable, then $X(T)\to\Pi_{X}(T)$ is an equivalence of categories.
\end{TC}

Recall that Grothendieck defined a geometrically connected variety as elementary anabelian if it can be constructed by successive smooth fibrations from hyperbolic curves, see \cite{gro97}. Merging the concepts of orbicurves and elementary anabelian varieties, in \autoref{elemsect} we define \emph{elementary anabelian stacks}. The section conjecture for curves implies that they are printable.

\begin{TD}[\ref{elem}]
	Elementary anabelian stacks over a field $k$ finitely generated over $\Q$ are fff. 
	
	If the section conjecture holds for smooth, proper, hyperbolic curves defined over fields finitely generated over $\Q$, then elementary anabelian stacks defined over fields finitely generated over $\Q$ are printable.
\end{TD}

Finally, we highlight two minor results that we think are worth observing.

\begin{itemize}
	\item If a smooth, proper DM stack $X/k$ is fff, then it has a finite étale cover by an algebraic space, see \autoref{cover}. This suggests that "anabelian" stacks should have a finite étale cover by an algebraic space.
	\item If a smooth, proper DM stack is printable, then $\pi_{1}(X_{\bar{k}})$ has no finite index abelian subgroups, see \autoref{anvsab}.
\end{itemize}

In \autoref{etaleapp} we prove some tools we need which are straightforward generalizations of the work of N. Borne and A. Vistoli in \cite{bv15}.

\subsection{Conventions and notations}

We always work over a field $k$ of characteristic $0$, except in \autoref{etaleapp} where there are no hypotheses on the base field. Curves and orbicurves will always be smooth, geometrically connected and proper, except if we specify differently.

We use underlines to distinguish between sets and sheaves: for instance, $\Pic$ is the Picard group, while $\underline{\Pic}$ is the Picard sheaf, or if $X$ is a stack with a rational point $x\in X(k)$, then $\uaut_{X}(x)$ is the sheaf of automorphisms of $x$, while $\aut_{X}(x)=\uaut_{X}(x)(k)$.

If $X$ is geometrically connected, we will denote by $\pi_{X}$ the structure morphism $X\to\Pi_{X/k}$ of the étale fundamental gerbe, see \cite{bv15} and \autoref{etaleapp}. If there is no risk of confusion, we may drop the subscript and write just $\pi:X\to\Pi_{X/k}$. We write $\pi_{1}(X,x)$ for classical étale fundamental groups and $\upi_{1}(X,x)=\uaut_{\Pi_{X/k}}(\pi(x))$ for étale fundamental group schemes.

Throughout the article (except in \autoref{open}), we restrict our attention to $X$ proper. The reason is that the section conjecture is much easier to handle in the proper case, and if one states the anabelian conjectures for DM stacks rather than schemes then the non-proper case can be recovered from the proper one using a limit process found by N. Borne and M. Emsalem (see \cite[§ 2.2.3]{be14} and \autoref{open}).

There is a small conflict of terminology between two of our major references. For N. Borne and A. Vistoli in \cite{bv15}, a finite stack over a field $k$ is a stack over $k$ which admits a presentation by a finite groupoid. A finite gerbe is a finite stack which is a gerbe. For the Stacks Project \cite{stacks-project}, finite morphisms are assumed to be representable. We stick with the Borne-Vistoli terminology.

\section{Stacky going up and going down theorems}

To understand precisely how anabelian geometry for DM stacks should look like, the single most important fact to understand is how the section conjecture behaves along finite étale morphism. In a classical context, i.e. for schemes, this situation is well understood and packed in the so called "going up" and "going down" theorems, see \cite[Propositions 110, 111]{sti13}. The formalism of étale fundamental gerbes is particularly well suited for the study of this situation: in fact, if $f:Y\to X$ is a representable, finite étale morphism, the natural diagram
\[\begin{tikzcd}
	Y\rar["\pi_{Y}"]\dar["f"]	&	\Pi_{Y/k}\dar["\pi_{f}"]	\\
	X\rar["\pi_{X}"]			&	\Pi_{X/k}
\end{tikzcd}\]
is $2$-cartesian, see \autoref{cartetale}. This fact makes the study of finite étale morphism with respect to the section conjecture particularly easy, even for stacks.

\begin{proposition}[Going up]\label{up}
    Let $X,Y$ be geometrically connected fibered categories and $f:Y\to X$ a representable, finite étale morphism. The following are true:
    \begin{enumerate}[(i)]
	    \item If $X(k)\to\Pi_{X/k}(k)$ is fully faithful, then $Y(k)\to\Pi_{Y/k}(k)$ is fully faithful, too.
	    \item If $X(k)\to\Pi_{X/k}(k)$ is an equivalence, then $Y(k)\to\Pi_{Y/k}(k)$ is an equivalence, too.
	\end{enumerate}
	\begin{proof}
		Follows directly from the fact that the diagram above is $2$-cartesian.
	\end{proof}
\end{proposition}

\begin{lemma}[Extension of the base field]\label{fieldext}
	Let $f:A\to B$ be a morphism of fibered categories over $k$ which are stacks in the étale topology, and $L/k$ a finite Galois extension. 
	\begin{enumerate}[(i)]
		\item Let $a,a'\in A(k)$ be rational points. If the map $\isom_{A}(a_{L},a'_{L})\to\isom_{B}(f(a_{L}),f(a'_{L}))$ is bijective, then the map $\isom_{A}(a,a')\to\isom_{B}(f(a),f(a'))$ is bijective.
	    \item If $A(L)\to B(L)$ is fully faithful, then $A(k)\to B(k)$ is fully faithful.
	    \item Let $b\in B(k)$ be a rational point, and suppose that $A(L)\to B(L)$ is fully faithful. Then $b$ is in the essential image of $A(k)\to B(k)$ if and only if $b_{L}$ is in the essential image of $A(L)\to B(L)$.
		\item If $A(L)\to B(L)$ is an equivalence, then $A(k)\to B(k)$ is an equivalence, too.
	\end{enumerate}
	\begin{proof}
		\begin{enumerate}[(i)]
	        \item We have a commutative diagram
	        \[\begin{tikzcd}
				\isom_{A}(a,a')\rar\dar[hook]				&	\isom_{B}(f(a),f(a'))\dar[hook]	\\
				\isom_{A}(a_{L},a_{L}')\rar["\sim"]			&	\isom_{B}(f(a_{L}),f(a_{L}))
	        \end{tikzcd}\]
	        where the vertical arrows are injective, and the lower arrow is bijective by hypothesis. Both $A$ and $B$ are stacks in the étale topology, hence the $\uisom$ functors are sheaves and satisfy Galois descent. This means that the sets in the upper row are just the $\gal(L/k)$-invariant elements of the groups in the lower row. Since the lower horizontal arrow is clearly equivariant, we get that the upper horizontal row is bijective, too.
			
			\item Follows from (i).
			
	        \item The "only if" part is obvious. Now suppose that $b_{L}\simeq f(a')$ is in the essential image of $A(L)\to B(L)$. For every $\sigma\in\gal(L/k)$, we have an isomorphism
	        \[\phi_{\sigma}:\sigma^{*}f(a')\simeq\sigma^{*}b_{L}=b_{L}\simeq f(a')\]
	        which corresponds to an isomorphism $\psi_{\sigma}:\sigma^{*}(a')\simeq a'$ since $A(L)\to B(L)$ is fully faithful by hypothesis.
	        
			Now, we have $\phi_{\sigma\rho}=\phi_{\sigma}\circ\sigma^{*}\phi_{\rho}$ by direct computation. Since $A(L)\to B(L)$ is fully faithful, this means that we also have $\psi_{\sigma\rho}=\psi_{\sigma}\circ\sigma^{*}\psi_{\rho}$ and hence by Galois descent there exists $a\in A(k)$ such that $a_{L}\simeq a'$. Let us check that $f(a)\simeq b$.
			
			We have a chain of isomorphisms
			\[f(a)_{L}=f(a_{L})\simeq f(a')\simeq b_{L},\]
			we have to check that this is Galois invariant. This amounts to the fact that, by definition, $f(\psi_{\sigma})=\phi_{\sigma}$.
			
	        \item Follows from (ii) and (iii).
		\end{enumerate}
	\end{proof}
\end{lemma}

In the following, we will use without mention the fact that, if $X$ is a geometrically connected fibered category and $L/k$ is a finite, separable extension, the natural morphism $\Pi_{X_{L}/L}\to\Pi_{X/k}\times_{k}L$ is an isomorphism (see \autoref{basechangefin}).

\begin{definition}
	Let $\mc{C},\mc{D}$ be categories and $f:\mc{C}\to \mc{D}$ a functor, $p\in \mc{C}$ an object. We say $f$ is fully faithful at $p$ if $\aut_{\mc{C}}(p)\to\aut_{\mc{D}}(f(p))$ is bijective.
\end{definition}

\begin{remark}
	Suppose that $\mc{C},\mc{D}$ are small categories in which all morphisms are isomorphisms. For example, $X(S)$ has this form for every stack $X$ and every scheme $S$. A functor $f:\mc{C}\to \mc{D}$ is fully faithful if and only if it is fully faithful at every point and is injective on isomorphism classes.
\end{remark}

\begin{proposition}[Going down]\label{down}
    Let $X$ and $Y$ be geometrically connected fibered categories which are stacks in the étale topology, and $f:Y\to X$ a representable, finite étale morphism. The following are true:
    \begin{enumerate}[(i)]
	    \item If $Y_{L}(L)\to\Pi_{Y_{L}}(L)$ is fully faithful for every finite, separable extension $L/k$, then $X(k)\to\Pi_{X/k}(k)$ is fully faithful.
	    \item If $Y_{L}(L)\to\Pi_{Y_{L}}(L)$ is an equivalence for every finite, separable extension $L/k$, then $X(k)\to\Pi_{X/k}(k)$ is an equivalence.
	\end{enumerate}
    \begin{proof}
		As in \autoref{up}, we are going to use the fact that the $2$-commutative diagram
		\[\begin{tikzcd}
			Y\rar["\pi_{Y}"]\dar["f"]	&	\Pi_{Y/k}\dar["\pi_{f}"]	\\
			X\rar["\pi_{X}"]			&	\Pi_{X/k}
		\end{tikzcd}\]
		is 2-cartesian, see \autoref{cartetale}. The proof is more complex than that of \autoref{up}, since now we have to make a descent argument.
        
        \begin{enumerate}[(i)]
			\item First, let us check that $X(k)\to\Pi_{X/k}(k)$ is fully faithful at every point, next we will show that it is injective on isomorphism classes.
	        
	        {\bf Full faithfulness at a point.} Choose $x\in X(k)$, since $Y\to X$ is finite étale there exists a finite Galois extension $L$ and a point $y\in Y_{L}(L)$ such that $f(y)\simeq x_{L}$. Thanks to \autoref{fieldext}.(i), we may suppose $L=k$, $f(y)\simeq x$. Write 
	        \[G=\aut_{\Pi_{X/k}}(\pi_{X}(x)),~~~H=\aut_{\Pi_{Y/k}}(\pi_{Y}(y)),\]
	        since $Y\to X$, $\Pi_{Y/k}\to\Pi_{X/k}$ are faithful we have natural embeddings
	        \[\pi_{f}:H\s G,~~~f:\aut_{Y}(y)\s\aut_{X}(x).\]
	        By an abuse of notation, write $\pi_{X}$, $\pi_{Y}$ for the homomorphisms $\aut_{X}(x)\to G$, $\aut_{Y}(y)\to H$ respectively.
	        
	        We have an isomorphism
	        \[\aut_{Y}(y)\simeq \aut_{X}(x)\times_{G}H,\]
	        and we also know that
	        \[\pi_{Y}:\aut_{Y}(y)\to H\]
			is an isomorphism. In particular, the fact that $\pi_{Y}:\aut_{Y}(y)\to H$ is injective implies that $\pi_{X}:\aut_{X}(x)\to G$ is injective. Let us prove surjectivity.
			
			Fix an element $g\in G$, since the diagram above is 2-cartesian the triple
			\[(x,\pi_{Y}(y),g)\]
			gives us a point $y'\in Y(k)$ such that $\pi_{Y}(y')\simeq\pi_{Y}(y)$ and $f(y')\simeq x$. Since $Y(k)\to\Pi_{Y/k}(k)$ is fully faithful, there exists an isomorphism $y\to y'$. Using the 2-cartesianity of the diagram above, the isomorphism $y\to y'$ gives us the following data: two isomorphisms $\alpha:x\to x$, $\alpha\in\aut_{X}(x)$ and $h:\pi_{Y}(y)\to\pi_{Y}(y)$, $h\in H$ such that
			\[\pi_{f}(h)\circ \id=\pi_{X}(\alpha)\circ g\in G.\]
			Since $\pi_{Y}:\aut_{Y}(y)\to H$ is an isomorphism, there exists $\beta\in\aut_{Y}(y)$ such that $\pi_{Y}(\beta)=h$. It follows that $g=\pi_{X}(\alpha^{-1}\circ f(\beta))$.
			
			{\bf Injectivity on isomorphism classes.} Suppose that we have an isomorphism $\alpha:\pi_{X}(x)\to\pi_{X}(x')$ for some $x,x'\in X(k)$, we want to show that there exists an isomorphism $x\to x'$, thanks to the preceding point this is equivalent to showing that 
			\[\isom_{X}(x,x')\to\isom_{\Pi_{X/k}}(\pi_{X}(x),\pi_{X}(x'))\]
			is bijective. There exists a finite Galois extension $L/k$ and a point $y\in Y(L)$ such that $f(y)=x_{L}$, thanks to \autoref{fieldext}.(i) we may assume $L=k$.
			
	        Since 
	        \[\pi_{f}(\pi_{Y}(y))=\pi_{X}(f(y))=\pi_{X}(x)\simeq\pi_{X}(x'),\]
	        by 2-cartesianity there exists a point $y'\in Y(k)$ such that $\pi_{Y}(y')\simeq\pi_{Y}(y)$ and $f(y')\simeq x'$. Now since $Y(k)\to\Pi_{Y/k}(k)$ is fully faithful by hypothesis and $\pi_{Y}(y)\simeq\pi_{Y}(y')$, we get an isomorphism $y\simeq y'$ which induces an isomorphism $x\simeq x'$ as desired.
	        
	        \item This is a direct application of point (i) and \autoref{fieldext}.(iii), together with the observation that every section $\spec k\to\Pi_{X/k}$ lifts to a section of $\Pi_{Y/k}$ up to a finite, separable field extension: in fact, $\spec k\times_{\Pi_{X/k}}\Pi_{Y/k}$ is a finite étale scheme. To check that $\spec k\times_{\Pi_{X/k}}\Pi_{Y/k}$ is finite étale, observe that up to an extension $k'/k$ we have
	        \[\spec k'\times_{\Pi_{X/k}}\Pi_{Y/k}\simeq \spec k'\times_{X}Y\]
	        for some point $\spec k'\to X$, since $\Pi_{X/k}$ is a gerbe and hence all points are fpqc locally isomorphic.
		\end{enumerate}
	\end{proof}
\end{proposition}

In the classical going up and down theorems there are hypotheses on the so called centralizers of sections. If $\sigma\in\Pi_{X/k}(k)$ corresponds to a section $s:\gal(\bar{k}/k)\to\pi_{1}(X,\bar{x})$, the centralizer of $s$ is the group of elements of $\pi_{1}(X_{\bar{k}},\bar{x})$ centralizing the image of $s$. However, in our results these hypotheses seem to be absent: the reason is that the notion of centralizer of a section (see \cite[§ 3.3]{sti13}) fits nicely in our point of view without any additional work. The following \autoref{centralizers} explains how.

\begin{lemma}\label{centralizers}
	Let $s:\gal(\bar{k}/k)\to\pi_1(X,x)$ be a section of the natural projection $\pi_1(X,x)\to\gal(\bar{k}/k)$, and $C_s\s\pi_1(X_{\bar{k}},x)$ its group of centralizers. Let $\sigma\in\Pi_{X/k}(k)$ the rational section corresponding to $s$. There is an isomorphism
	\[C_s\simeq\uaut_{\Pi_{X/k}}(\sigma)(k).\]
	\begin{proof}
		This follows from \autoref{fundsect}. Let us explain this.
		
		We have a natural identification
		\[\pi_{1}(X_{\bar{k}},\bar{x})=\uaut_{\Pi_{X/k}}(\pi(x))(\bar{k}).\]
		
		Since $\Pi_{X/k}$ is a gerbe, there exists an isomorphism $\Phi:\uaut_{\Pi_{X/k}}(\pi(x))(\bar{k})\simeq\uaut_{\Pi_{X/k}}(\sigma)(\bar{k})$. The section $s$ induces an action of $\gal(\bar{k}/k)$ on $\pi_{1}(X_{\bar{k}},x)$ by conjugation, and this action coincides with the natural action on $\uaut_{\Pi_{X/k}}(\sigma)(\bar{k})$ pulled back to $\uaut_{\Pi_{X/k}}(\pi(x))(\bar{k})=\pi_{1}(X_{\bar{k}},x)$. Hence $g\in\pi_{1}(X_{\bar{k}},x)$ centralizes $s$ if an only if $\Phi(g)\in\uaut_{\Pi_{X/k}}(\sigma)(\bar{k})$ is Galois invariant, i.e. it is rational.
	\end{proof}
\end{lemma}

\section{Printable DM stacks}\label{s2}

Now that we have established what happens along finite, étale covers, we want to understand what the section conjecture for DM stacks should look like. Clearly, one can just directly translate Grothendieck's section conjecture to DM stacks. Here we hope to show that the right thing to conjecture in general is slightly stronger (but equivalent in the case of hyperbolic curves).

\begin{proposition}\label{strong}
	Let $X$ be a proper, smooth, geometrically connected Deligne-Mumford stack over $k$. The following are equivalent:
	\begin{enumerate}
		\item for every finitely generated extension $k'/k$ and for every finite étale geometrically connected cover $Y\to X_{k'}$,
		\[Y(k')\to\on{Hom-ext}_{\Gamma_{k'}}(\Gamma_{k'},\pi_{1}(Y))\]
		is bijective (resp. injective) on isomorphism classes,
		\item the natural map
		\[X(k')\to\Pi_{X/k}(k')\]
		is an equivalence of categories (resp. fully faithful) for every finitely generated extension $k'/k$.
	\end{enumerate}
	\begin{proof}
		Suppose that $X(k')\to\Pi_{X/k}(k')$ is an equivalence (resp. fully faithful). Then by \autoref{basechange} $X_{k'}(k')\to\Pi_{X_{k'}/k'}(k')$ is an equivalence (resp. fully faithful), too, and hence $Y(k')\to\on{Hom-ext}_{\Gamma_{k'}}(\Gamma_{k'},\pi_{1}(Y))$ is bijective (resp. injective) thanks to the going up theorem \autoref{up}.
		
		Suppose now that (1) holds, let $k'/k$ be a finitely generated extension, $x\in X(k')$ a point and $\pi(x)\in\Pi_{X/k}(k')$. Write $G=\uaut_{\Pi_{X/k}}(\pi(x))$. Since by hypothesis $X(k')\to\Pi_{X/k}(k')$ is bijective (resp. injective) on isomorphism classes, then we only have to show that
		\[\uaut_{X}(x)\to G=\uaut_{\Pi_{X/k}}(\pi(x))\]
		induces a bijection on $k'$-rational points. Thanks to \autoref{basechange}, we may suppose $k'=k$.
		
		{\bf Surjectivity of $\uaut_{X}(x)(k)\to G(k)$.} Let $A\s G(k)$ be the image of $\uaut_{X}(x)(k)$, and assume by contradiction that $g\notin A$. Since $A\s G(k)$ is finite and $g\notin A$, there exists a finite index subgroup $H\s G$ such that $A\s H(k)$ and $g\notin H(k)$: in order to find it, choose a finite quotient $q:G\to Q$ such that $q(g)\notin q(A)$, then choose $H$ as the inverse image of $q(A)$. Now consider the 2-fiber product
		\[\begin{tikzcd}
			Y\rar\dar	&	BH=\Pi_{Y/k}\dar	\\
			X\rar		&	BG=\Pi_{X/k}
		\end{tikzcd}\]
		where $BH$ identifies naturally with $\Pi_{Y/k}$. In fact, the universal property of $\Pi_{Y/k}$ gives us a natural map $\Pi_{Y/k}\to BH$, and thanks to \autoref{cartetale} $\Pi_{Y/k},BH$ are both subgerbes of $\Pi_{X/k}$ with the same finite index, hence they coincide.
		
		Now consider $x\in X(k)$ and the distinguished point $d_{H}:\spec k\to BH$: by construction, they both map to $\pi(x)\in\Pi_{X/k}(k)$, i.e. the distinguished point of $BG$. By definition of 2-fiber product, every automorphism of $\pi(x)$ defines a rational point of $Y$ mapping to $x\in X(k)$ and $d_{H}\in BH(k)$. In particular, we have two such rational points $y=(x,d_{H},\id)$, $y'=(x,d_{H},g)$ in $Y(k)$, let us study the isomorphisms $y\to y'$.
		
		By definition of 2-fiber product, an isomorphism $\beta:y=(x,d_{H},\id)\to y'=(x,d_{H},g)$ is given by a couple of isomorphisms $\alpha:x\to x$, $\alpha\in\uaut_{X}(x)$, and $h:d_{H}\to d_{H}$, $h\in H(k)$, such that 
		\[h\circ \id=g\circ\pi(\alpha)\in G(k).\]
		By construction, $\pi(\alpha)\in A\s H(k)$ and clearly $h\in H(k)$: since we have chosen $H$ such that $g\notin H(k)$, the equation above tells us that an isomorphism $y\to y'$ cannot exists, i.e. $y,y'$ are not isomorphic. But this gives an absurd, since $y,y'$ both map to $d_{H}\in BH(k)$ and $Y(k)\to BH(k)=\Pi_{Y/k}(k)$ is injective on isomorphism classes by hypothesis.
		
		{\bf Injectivity of $\uaut_{X}(x)(k)\to G(k)$.} Since $\uaut_{X}(x)$ is finite étale, up to enlarging the base field we may suppose that $\uaut_{X}(x)$ is discrete. Let $A\s G(k)$ be the image of $\uaut_{X}(x)(k)$. We can find a finite index subgroup $H\s G$ such that
		\[H(k)\cap A=\{\id\}\s G(k),\]
		for instance by taking a finite quotient $G\to Q$ such that $A\to Q(k)$ is injective and choosing $H$ as the kernel.
		
		Take the 2-fiber product $Y=X\times_{\Pi_{X/k}}BH$ as above, $x\in X(k)$ and the distinguished point $d_{H}:\spec k\to BH$ define a rational point $y=(x,d_{H},\id)\in Y(k)$. Since
		\[\uaut_{Y}(y)\simeq \uaut_{X}(x)\times_{G}H\]
		and $A\cap H(k)=\{\id\}$, we get that $\uaut_{Y}(y)\s\uaut_{X}(x)$ is the kernel of $\pi_{X}:\uaut_{X}(x)\to G$. Hence, we want to prove that $\uaut_{Y}(y)$ is trivial.
		
		Suppose by contradiction that $\uaut_{Y}(y)$ is not trivial. A non-trivial finite étale group scheme is non-special, hence there exists a finitely generated extension $k'/k$ and a point $y':\spec k'\to B\uaut_{Y}(y)$ which is not $k'$-isomorphic to the distinguished one, i.e. $y'\not\simeq y\in B\uaut_{Y}(y)(k')$.
	
		The morphism $B\uaut_{Y}(y)\to Y$ is the residual gerbe (see \cite[\href{https://stacks.math.columbia.edu/tag/06MU}{Definition 06MU}]{stacks-project}) of $y$ thanks to \cite[\href{https://stacks.math.columbia.edu/tag/0DTI}{Lemma 0DTI}, \href{https://stacks.math.columbia.edu/tag/06UI}{Lemma 06UI}]{stacks-project}, in particular it is a monomorphism and thus $y,y'$ define non-isomorphic points of $Y(k')$. On the other hand, since
		\[\uaut_{Y}(y)\to H=\uaut_{\Pi_{Y/k}}(\pi_{Y}(y))\s G\]
		is a trivial homomorphism of group schemes, the images of $y,y'$ in $BH(k')=\Pi_{Y/k}(k')$ are isomorphic. This gives an absurd, since by hypothesis $Y(k')\to\Pi_{Y/k}(k')$ is injective on isomorphism classes.
	\end{proof}
\end{proposition}

We define \emph{printable DM stacks} as those ones satisfying the equivalent conditions of \autoref{strong}.

\begin{definition}\label{strongdef}
	Let $X$ be a smooth, proper, geometrically connected Deligne-Mum\-ford stack. We say that $X$ is \emph{printable} (resp. \emph{fundamentally fully faithful, or fff}) if the natural morphism
	\[X(k')\to\Pi_{X/k}(k')\]
	is an equivalence of categories (resp. fully faithful) for every finitely generated extension $k'/k$. 
\end{definition}

As we will see later, even if this definition seems deeply arithmetic in nature, it is actually purely geometric: if $k\s\C$, printability of $X$ depends only on $X_{\C}$, see \autoref{anabgeom}. This agrees with the ideas expressed by Grothendieck in \cite{gro97}.

\begin{remark}
	Extending the definition to Deligne-Mumford stacks seems natural for at least two reasons. One is that moduli stacks of curves are expected to be anabelian, the second is that hyperbolic orbicurves are printable if and only if hyperbolic curves are printable, see \autoref{eq}. We address the question "why not Artin stacks?" in \autoref{artinsect}.
\end{remark}

Let us study printability in dimension $0$.

\begin{lemma}\label{redgerbe}
	A geometrically connected, geometrically reduced $0$-dimensional DM stack of finite type over $k$ is a finite étale gerbe.
	\begin{proof}
		Let $X$ be such a DM stack. Applying the definition of gerbe, it is immediate to see that $X$ is a gerbe over $k$ if and only if $X_{\bar{k}}$ is a gerbe over $\bar{k}$, hence we may assume $k$ algebraically closed. Let $U\to X$ be an étale cover of finite type and $R=U\times_{X}U$, we have that $R\rightrightarrows U$ is a groupoid and the natural map $[U/R]\to X$ is an isomorphism.
		
		Since $U,R$ are $0$ dimensional, reduced schemes of finite type over the algebraically closed field $k$, they are finite disjoint unions of copies of $\spec k$, i.e. we may think of them as sets. The groupoid $R\rightrightarrows U$ induces an oriented graph whose set of nodes is $U$ and whose set of arrows is $R$, let $(U,R)=(U_{1},R_{1})\sqcup\dots\sqcup (U_{n},R_{n})$, be the connected components of the graph. For every $i$, we have an induced groupoid $R_{i}\rightrightarrows U_{i}$, and by construction $X=\bigsqcup_{i}[U_{i}/R_{i}]$. Since $X$ is connected, it follows that $n=1$.
		
		Now, since the graph induced by $U,R$ is connected, it is immediate to check that $X=[U/R]=BG$ where $G$ is the group of $R$-automorphisms of any point of $U$.
	\end{proof}
\end{lemma}

\begin{corollary}
	A smooth, proper, geometrically connected DM stack of dimension $0$ over $k$ is printable.
	\begin{proof}
		If $X$ is such a DM stack, it is a finite étale gerbe thanks to \autoref{redgerbe}, thus $X=\Pi_{X/k}$ is obviously printable.
	\end{proof}
\end{corollary}

In the following, we show what it means for a scheme to be printable in the classical terms of the section conjecture and of centralizers of sections, see \cite[§ 3.3]{sti13}.

\begin{lemma}\label{schprint}
	Let $X$ be a smooth, proper, geometrically connected scheme. Then $X$ is printable (resp. fff) if and only if
	\begin{itemize}
		\item $X_{k'}$ satisfies the section conjecture (resp. the injectivity part of the section conjecture) for every finitely generated extension $k'/k$, and
		\item for every $x\in X(k')$, the associated section in $\on{Hom-ext}_{\Gamma_{k'}}(\Gamma_{k'},\pi_{1}(X))$ has trivial centralizer.
	\end{itemize}
	\begin{proof}
		As we have shown in \autoref{centralizers}, the automorphism groups of the points of the fundamental gerbe correspond to centralizers of sections of the étale fundamental group, hence if $X$ is a scheme asking an equivalence of categories corresponds to asking a bijection on isomorphism classes together with the triviality of centralizers.
	\end{proof}
\end{lemma}

\begin{proposition}\label{curvessect}
	Let $k$ be finitely generated over $\Q$.
	\begin{itemize}
		\item Smooth proper curves over $k$ are fundamentally fully faithful if and only if they have positive genus.
		\item Hyperbolic curves over $k$ are printable if and only if they satisfy the section conjecture over every finitely generated extension of the base field.
	\end{itemize}
	\begin{proof}
		For smooth, proper curves with Euler characteristic less than or equal to $0$, centralizers of sections coming from rational points are trivial, thanks either to \cite[Proposition 36, Proposition 104]{sti13} or to the full faithfulness part of \autoref{strong}. Apply \autoref{schprint}.
	\end{proof}
\end{proposition}

\begin{proposition}\label{etale}
	Let $Y$, $X$ be smooth, proper, geometrically connected DM stacks over a field $k$, and $Y\to X$ a finite étale covering. Then $X(k')\to\Pi_X(k')$ is an equivalence (resp. fully faithful) for every finite, separable  extension $k'/k$ if and only if the same holds for $Y$. In particular, $Y$ is printable (resp. fff) if and only if $X$ is printable (resp. fff).
	\begin{proof}
		This is a straightforward application of the going up and down theorems \autoref{up}, \autoref{down}.
	\end{proof}
\end{proposition}

\section{Why not Artin stacks}\label{artinsect}

One may wonder: why DM stacks and not algebraic (i.e. Artin) stacks? The answer is based on one's taste. DM stacks seem more natural, since $\Pi_{X/k}$ is profinite étale and \autoref{strong} fails for algebraic stacks. For example, if $G$ is a connected algebraic group, then condition (1) of \autoref{strong} holds for $BG$ if an only if $G$ is special, while condition (2) if and only if $G$ is trivial. Hence, it makes a difference if we choose condition (1) or (2) as definition of printability for algebraic stacks.

If we choose (1), we should for instance consider $B\on{GL}_{n}$ as printable even if $B\on{GL}_{n}\to\Pi_{B\on{GL}_{n}}=\spec k$ is not an equivalence of categories on rational points, and this seems not very pleasant. On the other hand, if we choose (2), the following proposition shows that we get back to DM stack.

\begin{proposition}\label{artin}
	Let $X$ be a separated, geometrically connected algebraic stack locally of finite type over $k$. Suppose that
	\[X(k')\to\Pi_{X/k}(k')\]
	is fully faithful for every finitely generated extension $k'/k$. Then $X$ is a DM stack.
	\begin{proof}
		Since we are in characteristic $0$, it is enough to show that $\uaut_{X}(x)$ is finite for any geometric point $x$, see \cite[\href{https://stacks.math.columbia.edu/tag/0DSN}{Lemma 0DSN}]{stacks-project}. Since $X$ is locally of finite type, we may assume that $x$ is defined over a finitely generated extension $k'/k$. Thanks to \autoref{basechange}, we may suppose $k'=k$, i.e. $x\in X(k)$ is a rational point. Since $X$ is separated, $\uaut_{X}(x)$ is a group scheme of finite type, see \cite[\href{https://stacks.math.columbia.edu/tag/0DTS}{Lemma 0DTS}]{stacks-project}.
		
		Let $\pi(x)\in\Pi_{X/k}(k)$ be the image of $x$, we have an homomorphism of group schemes
		\[\uaut_{X}(x)\xar{\pi}\uaut_{\Pi_{X/k}}(\pi(x)).\]
		This homomorphism has trivial kernel: otherwise, since $\uaut_{X}(x)$ is of finite type, up to enlarging the base field we may suppose that there exists a rational point $\phi\in\ker(\pi)(k)$ different from the identity. But $\uaut_{X}(x)\to\uaut_{\Pi_{X/k}}(\pi(x))$ is injective on rational points by hypothesis, we are in characteristic $0$ and hence $\ker(\pi)$ is trivial. By the following \autoref{fpf}, we get that $\uaut_{X}(x)$ is finite.
	\end{proof}
\end{proposition}

\begin{lemma}\label{fpf}
	Let $f:G\to P$ be an homomorphism of group schemes over a field $k$ with trivial kernel. Assume that $G$ is of finite type and $P$ is pro-finite. Then $G$ is finite.
	\begin{proof}
		We may base change everything to $\bar{k}$ and assume that $k$ is algebraically closed. Let $g\in G$ be any point with residue field $k(g)$. Since $G$ is of finite type, there exists an irreducible variety $U$ with function field $k(U)=k(g)$ and a locally closed embedding $U\s G$. Since $P$ is pro-finite and $k$ is algebraically closed, $f(g)$ is rational, thus by construction all the points of $U(k)\s G(k)$ map to $f(g)$. Since $G\to P$ has trivial kernel, it follows that $U(k)$ has only one point and hence $k(g)=k(U)=k$. Since $G$ is of finite type and every point is rational, it is finite.
	\end{proof} 
\end{lemma}

\section{Covers by algebraic spaces}\label{coversect}

It turns out that fff DM stacks must have a non-obvious topological feature: they are uniformizable in the sense of Noohi, i.e. they have a finite étale cover by an algebraic space, see \cite[Definition 6.1]{noo04}. Noohi essentially proves this in \cite[Theorem 6.2]{noo04}, but the connection with the section conjecture is not stated in his work.

Rather than using Noohi's result, we prove it again in our setting: a formal comparison with Noohi's theory would be longer. The idea behind the proof is essentially the same.

\begin{proposition}\label{cover}
	Let $X$ be a geometrically connected, separated DM stack of finite type over $k$, and suppose that the natural morphism
	\[X(k')\to\Pi_{X/k}(k')\]
	is faithful for every finitely generated $k'/k$.
	
	There exists a finite étale gerbe $\Phi$ with a representable morphism $X\to\Phi$, and a finite étale cover $E\to X$ with $E$ an algebraic space. If $\Pi_{X/k}(k)\neq\emptyset$, we can choose $E$ to be geometrically connected. 
	\begin{proof}
		Choose $\xi:\spec k'\to X$ the generic point of an irreducible component of $X$ with $k'$ finitely generated over $k$. Up to enlarging $k'$, we may assume that $\uaut_{X}(\xi)$ is discrete over $k'$. By hypothesis, there exists a finite étale gerbe $\Phi_{0}$ over $k$ with a morphism $\phi_{0}:X\to\Phi_{0}$ which is faithful at $\xi$, i.e. $\uaut_{X}(\xi)\to\uaut_{\Phi_{0}}(\phi_{0}(\xi))$ is injective.
		
		Consider the relative inertia $I_{X/\Phi_{0}}$. By generic flatness, there exists an open, irreducible subset $U_{0}\s X$ such that $\xi\in U_{0}$ and the restriction of $I_{X/\Phi_{0}}\to X$ to $U_{0,\rm red}\s X_{\rm red}$ is flat. Since $I_{X/\Phi_{0}}\to X$ is proper and unramified, too, its restriction to $U_{0,\rm red}$ is finite étale. Since $X\to \Phi_{0}$ is faithful at $\xi$, it moreover follows that the restriction of $I_{X/\Phi_{0}}$ to $U_{0,\rm red}$ is a finite étale morphism of degree $1$, i.e. an isomorphism. In particular, $I_{U_{0}/\Phi_{0}}\to U_{0}$ is radicial and unramified, thus a monomorphism. Moreover, it has a section, hence it is an isomorphism and hence $U_{0}\to\Phi_{0}$ is faithful.
		
		Now let $X_{1}=X\setminus U_{0}$, we may repeat the process and find $U_{1}\s X_{1}$ with a finite étale gerbe $\Phi_{1}$ and a morphism $X\to\Phi_{1}\to\Phi_{0}$ such that $U_{1}\to \Phi_{1}$ is faithful, define $X_{2}=X_{1}\setminus U_{1}$ and continue by recursion. By noetherian descent, the process ends, thus for some large $N$ we have that $X_{N+1}$ is empty and $X\to\Phi_{N}$ is faithful. Choose $\Phi=\Phi_{N}$. We have that $X\to\Phi$ is representable since a faithful morphism of algebraic stacks is representable, see \cite[\href{https://stacks.math.columbia.edu/tag/04Y5}{Lemma 04Y5}]{stacks-project}.

		In order to find $E$, since $\Phi$ is a finite étale gerbe there exists a finite, separable extension $k'/k$ and a section $\spec k'\to\Phi$. Take $E=\spec k'\times_{\Phi}X$.
		
		Now suppose that $\Pi_{X/k}(k)\neq\emptyset$, in particular we have a section $\spec k\to\Phi$. Thanks to \cite[Lemma 5.12]{bv15}, we may assume that $X\to\Phi$ is Nori-reduced: this exactly means that $E=\spec k\times_{\Phi}X$ is geometrically connected, see \cite[Remark 5.11]{bv15}. We remark that in \cite[Remark 5.11]{bv15} it is assumed that $X$ is geometrically reduced, an hypothesis we do not have, but it can be checked that they do not actually use it (provided that in their proof algebraic closures are replaced by separable closures).
	\end{proof}
\end{proposition}

\section{Printability depends only on the geometric type}

A priori, our definition of printable DM stack depends on the base field $k$. It turns out that it is actually independent of it.

In order to prove this, we have to generalize some Galois-theoretic facts about finite extensions of fields to finitely generated extensions. We work in characteristic $0$ in order to avoid inseparability issues.

\begin{lemma}\label{fggalext}
	Let $k'/h/k$ be finitely generated extensions of a field $k$ of characteristic $0$. Let $\sigma$ be a non-trivial automorphism of $h$ of finite order which is trivial on $k$. Then there exists a finite extension $k''/k'$ such that $\sigma$ extends to an automorphism of finite order of $k''$.
	\begin{proof}
		Up to replacing $k$ with the subfield of $h$ fixed by $\sigma$, we may assume that $h/k$ is finite.
	
		Let $t_{1},\dots,t_{n}$ be a transcendence basis of $k'/h$, we have that $\sigma$ extends to an automorphism of $h(t_{1},\dots,t_{n})$ trivial on $k(t_{1},\dots,t_{n})$ by acting trivially on $t_{i}$. Since $h$ is finite over $k$, the extension $k'/k(t_{1},\dots,t_{n})$ is finite. Choose a Galois closure $k''/k'/k(t_{1},\dots,t_{n})$, the automorphism of $h(t_{1},\dots,t_{n})$ extends to $k''$.
	\end{proof}
\end{lemma}

\begin{lemma}\label{fggalois}
	Let $k'/h/k$ be finitely generated extensions of a field $k$ of characteristic $0$. If $h/k$ is non-trivial, then there exists a finitely generated extension $k''/k'/k$ with an automorphism $\sigma:k''\to k''$ of finite order which is trivial on $k$ but non-trivial on $h$.
	\begin{proof}
		Let us do this in three cases.
		
		{\bf Case 1.} $k'=h$ and $k$ is not algebraically closed in $h$. Let $l=\bar{k}^{h}$, fix $\tilde{l}$ a Galois closure of $l$ and call $\tilde{h}=\tilde{l}\otimes_{l}h$. Since $l$ is algebraically closed in $h$ we get that $\tilde{h}$ is a domain, moreover it is a field since it is finite over $h$. Now take any automorphism $\sigma$ of $\tilde{l}/k$ which is non-trivial on $h\cap l$ (it exists by Galois theory since the extension $h\cap l/k$ is non-trivial) and apply \autoref{fggalext}.
		
		{\bf Case 2.} $k'=h$ and $k$ is algebraically closed in $h$. Let $k''$ be the fraction field of the domain $h\otimes h$, and let $\sigma$ be the automorphism of $h\otimes h$ which permutes the two coordinates: this extends to the fraction field $k''$, and thus we conclude.
		
		{\bf General case.} By the preceding cases, there exists a finitely generated extension $h'/h$ with an automorphism $\sigma$ of $h'$ trivial on $k$ but not trivial on $h$. Choose any finitely generated extension $k''/k$ which contains both $h'$ and $k'$ as subextensions. Thanks to \autoref{fggalext}, there exists a further finite extension $k'''/k''$ such that $\sigma$ extends to $k'''$.
	\end{proof} 
\end{lemma}

\begin{lemma}\label{rigid}
	Let $G$ be a profinite étale group scheme over $k$, and suppose that
	\[G(k')=\{\id\}\]
	for every field $k'$ finite over $k$. Let $T\to\spec k$ be a $G$-torsor and $k'/k$ a finitely generated extension such that $T_{k'}\to\spec k'$ is the trivial torsor. Then $T$ is trivial.
	\begin{proof}
		Let $p\in T$ be the image of a point $\spec k'\to T$, $k(p)/k$ is finite and separable since $k'/k$ is finitely generated and $T\to\spec k$ is profinite étale. Let $\tilde{k(p)}$ be a Galois closure of $k(p)/k$. We have that $\spec k(p)\otimes_k\tilde{k(p)}\to T_{\tilde{k(p)}}$ is a closed embedding because it is the base change of $\spec k(p)\to T$ which is a closed embedding. But if $k(p)/k$ is non-trivial $\spec k(p)\otimes_k\tilde{k(p)}$ is a finite étale scheme with more than one point, hence we get an absurd because $T_{\tilde{k(p)}}\simeq G_{\tilde{k(p)}}$ has only one rational point by hypothesis.
	\end{proof}
\end{lemma}

\begin{example}
	Let $X/k$ be an fff algebraic space, and $p\in X(k)$ a rational point. Then
	\[\upi_1(X,x)=\uaut_{\Pi_{X/k}}(\pi(x))\]
	respects the hypothesis of \autoref{rigid} by definition of fff. In classical terms, the fact that $\upi_1(X,x)$ has no non-trivial rational points amounts to the triviality of centralizers of $\pi(x)$, see \autoref{centralizers}.
\end{example}

\begin{lemma}\label{fgext}
	Let $k'/k$ be a finitely generated extension of the base field, and $X$ a separated, fff DM stack of finite type over $k$. Let $s\in\Pi_{X/k}(k)$ be such that $s'=s_{k'}$ is in the essential image of $X(k')\to\Pi_{X/k}(k')$. Then $s$ is in the essential image of $X(k)\to\Pi_{X/k}(k)$.
	\begin{proof}
		If $k'/k$ is finite, since we are in characteristic $0$ we may extend $k'/k$ and suppose it is Galois, then this is the content of \autoref{fieldext}.iii, hence we may freely take finite extensions of the base field.
		
		Thanks to \autoref{cover}, we may take a finite cover $E\to X$ with $E$ an algebraic space. Up to a finite extension of the base field and by taking connected components, we may suppose that $E$ is geometrically connected. Since every point of $X$ extends to $E$ up to a finite extension, by an easy diagram chasing we may replace $X$ with $E$ and suppose that $X$ is an algebraic space.
		
		Now, by hypothesis we have a point $x'\in X(k')$ which maps to $s'\in\Pi_{X/k}(k')$. Since $X$ is separated, $x'$ has a residue field $k(x')$. If the extension $k(x')/k$ is non-trivial, then thanks to \autoref{fggalois} there exists a finitely generated extension $k''/k'/k(x')/k$ with a non-trivial automorphism $\sigma:k''\to k''$ which fixes $k$ but does not fix $k(x')$. In particular, $\sigma^* x'_{k''}\neq x'_{k''}$, but
		\[\pi_{X}(\sigma^*x'_{k''})=\sigma^*\pi_{X}(x'_{k''})=\sigma^*s'_{k''}=s'_{k''}=\pi_{X}(x')\in\Pi_{X/k}(k'')\]
		because $s'$ is defined on $k$. This gives an absurd since $X$ is fff by hypothesis, hence $k(x')=k$ i.e. there exists $x\in X(k)$ with $x_{k'}=x'$. 
		
		We want now to show that $\pi_{X}(x)=s$ using the fact that 
		\[\pi_{X}(x)_{k'}=\pi_{X}(x_{k'})=\pi_{X}(x')=s_{k'}.\]
		We may think $s$ as a $\upi_{1}(X,x)=\uaut_{\Pi_{X/k}}(x)$-torsor, $\pi_{X}(x)$ is the trivial $\upi_{1}(X,x)$-torsor and $k'$ splits the torsor $s$. Then $s$ is trivial because $\upi_{1}(X,x)$ respects the hypothesis of \autoref{rigid} since $X$ is fff.
	\end{proof}
\end{lemma}

\begin{theorem}\label{anabext}
	Let $k'/k$ be a finitely generated extension, and $X$ a smooth, proper, geometrically connected DM stack over $k$. Then $X$ is printable (resp. fff) if and only if $X_{k'}$ is printable (resp. fff).
	
	As a consequence, if $k\s\C$ is finitely generated over $\Q$, whether or not $X$ is printable depends only on $X_{\C}$.
	\begin{proof}
		We only do this for printability, the argument for fff is strictly contained.
		
		If $X$ is printable, $X_{k'}$ is printable by definition since $\Pi_{X_{k'}/k'}=\Pi_{X/k}\times_{k}k'$ thanks to \autoref{basechange}.
		
		On the other hand, suppose that $X_{k'}$ is printable. If $k'/k$ is finite, up to a finite extension we may suppose that it is Galois, too. Then this is the content of \autoref{fieldext}.
		
		Now that we can take arbitrary finite extensions of the base field, we may reduce to the case in which $X$ is an algebraic space using the same argument we have used in \autoref{fgext}.
		
		Hence, we may suppose that $X$ is an algebraic space. Let $L/k$ be a finitely generated extension, we want to show that $X(L)\to\Pi_{X/k}(L)$ is an equivalence. There exists a finitely generated extension $L'$ of $k'$ containing $L$, up to extensions we may suppose $L=k$ and $L'=k'$.
		
		First, we must show that $\pi_{X}:X(k)\to\Pi_{X/k}(k)$ is fully faithful. Since $X$ is an algebraic space, this amounts to showing injectivity on isomorphism classes together with the fact that for every $x\in X(k)$, $\uaut_{\Pi_{X/k}}(\pi_{X}(x))(k)$ is trivial. But these are direct consequences of the analogous facts over $k'$, which are true by hypothesis.
		
		Finally, we have to show essential surjectivity of $\pi_{X}:X(k)\to\Pi_{X/k}(k)$, and this is the content of \autoref{fgext}.
		
		If $k\s\C$ is finitely generated over $\Q$, and $X,X'$ are DM stacks over $k$ with an isomorphism $X_{\C}\simeq X'_{\C}$, the isomorphism descends to a finitely generated extension of $k$, thus $X$ is printable if and only if $X'$ is printable.
	\end{proof}
\end{theorem}

\begin{remark}\label{anabgeom}
	Thanks to \autoref{anabext}, we can see printability as a geometric property, rather than an arithmetic one, and this is coherent with Grothendieck's ideas. Clearly this is a tautology, we are not really able to describe in purely geometrical terms which DM stacks over $\C$ descend to printable DM stacks: still, we think it is worth observing that the arithmetic property depends only on the geometry of the variety.
\end{remark}

\section{Orbicurves}\label{orbisect}

The first non-trivial example of expected anabelian DM stacks are hyperbolic orbicurves. Borne and Emsalem conjectured \cite[Conjecture 2]{be14} that the section conjecture holds for them.

A proper orbicurve is a smooth, proper DM stack of dimension $1$ which is generically a scheme. In order to be clear, we will use Fraktur letters for orbicurves and normal ones for schemes.

Let $\mf{X}$ be an orbicurve with coarse moduli space $\mf{X}\to X$, there exists a maximal open subset $U\s\mf{X}$ which is a scheme, $U\to X$ is an open immersion. Let $D=X\setminus U$, define the rational Euler characteristic of $\mf{X}$ as
\[\rchi(\mf{X})=2-2g-\sum_{x\in D}\frac{r_{x}-1}{r_{x}}[k(x):k]\]
where $r_{x}$ is the degree of the residual gerbe at $x$, i.e. the ramification degree of $\mf{X}\to X$ at $x$. The orbicurve $\mf{X}$ is hyperbolic (resp. elliptic, parabolic) if $\rchi(\mf{X})<0$ (resp. $=0$, $>0$), see \cite[§ 2.2]{be14}, \cite[Proposition 5.11]{bn06}.

If $\mf{Y}\to\mf{X}$ is a finite étale cover of degree $d$, Riemann-Hurwitz applied to the associated (possibly ramified) morphism of coarse moduli spaces $Y\to X$ yields to the usual formula
\[\rchi(\mf{Y})=d\rchi(\mf{X}).\]

\begin{proposition}\label{orbicover}
	Let $\mf{X}$ be an orbicurve with coarse moduli space $\mf{X}\to X$, $U\s\mf{X}$ maximal open subset which is scheme, $D=X\setminus U$, $r_{x}$ the degree of the residual gerbe at $x$ for $x\in D$. Suppose that we are not in one of the following cases:
	\begin{itemize}
		\item $g(X)=0$, $\deg D=1$;
		\item $g(X)=0$, $D=\{x_{1},x_{2}\}$, $x_{i}\in X(k)$ rational, $r_{x_{1}}\neq r_{x_{2}}$.
	\end{itemize}
	Then there exists a finite extension $k'/k$ and a smooth, geometrically connected curve $Y$ over $k'$ with a geometrically Galois finite étale cover $Y\to\mf{X}_{k'}$. If $\Pi_{\mf{X}/k}(k)\neq\emptyset$, we may furthermore assume $k'=k$.
	\begin{proof}
		Since everything is of finite type, with standard arguments we can obtain the general case once we know the proposition is true for $k$ finitely generated over $\Q$. Now we may fix an embedding $k\s\C$, and thus reduce to $k=\C$. For $k=\C$, this is \cite[Proposition 5.7]{bn06}.
	\end{proof}
\end{proposition}

We are able to work with orbicurves because those ones with non-positive Euler characteristic have a finite étale covering which is a curve, i.e. they are uniformizable in the sense of Noohi. This is not only an useful feature, but a necessary one: thanks to \autoref{strong} and \autoref{cover}, it is implied by the injectivity part of the section conjecture for orbicurves. It is rather remarkable that this necessary topological feature happens to be true.

\begin{theorem}\label{eq}
	Let $k$ be finitely generated over $\Q$.
	\begin{itemize}
		\item A smooth, proper orbicurve $\mf{X}$ is fundamentally fully faithful if and only if $\rchi(\mf{X})\le 0$.
		\item If smooth, proper, hyperbolic curves satisfy the section conjecture over every finite extension $k'/k$, the section map is an equivalence for smooth, proper, hyperbolic orbicurves over $k$.
		\item In particular, smooth, proper, hyperbolic orbicurves are printable if and only smooth, proper, hyperbolic curves are printable.  
	\end{itemize}
	\begin{proof}
		Thanks to \autoref{curvessect}, \autoref{etale}, \autoref{anabext} and \autoref{orbicover}, we may reduce to one of the following cases: $\mf{X}$ is either a curve or a simply connected orbicurve. Both these cases are obvious.
	\end{proof}
\end{theorem}

\begin{corollary}\label{classeq}
	The section conjecture holds for hyperbolic orbicurves if and only if it holds for hyperbolic curves.\qed
\end{corollary}

\section{Affine curves}\label{open}

There is a version of the section conjecture for affine curves. If $U$ is a smooth geometrically connected curve with smooth completion $X$ and complement $D=X\setminus U$, every "missing" rational point $x\in D(k)$ defines a so called \emph{packet} of cuspidal sections $\mc{P}_{x}\s\Pi_{U/k}(k)$, see \cite{eh08} and \cite{sti12}. Following \cite{sti12}, the packet based at the cusp $x$ is defined as follows. Let $\O^{\rm h}_{X,x}$ be the henselianization of the local ring, and define the scheme of nearby points $U^{x}=U\times_{X}\spec(\O^{\rm h}_{X,x})$. The packet $\mc{P}_{x}$ is then defined as the image of
\[\Pi_{U^{x}/k}(k)\to\Pi_{X/k}(k).\]

The section conjecture for $U$ says that if $k$ is finitely generated over $\Q$ and $U$ has negative Euler characteristic, then the section map
\[U(k)\sqcup \bigsqcup_{x\in D(k)}\Pi_{U^{x}/k}(k)\to\Pi_{U/k}(k)\]
is bijective on isomorphism classes.

As showed by Niels Borne and Michel Emsalem in \cite[§ 2.2.3]{be14}, the section conjecture for orbicurves implies easily the section conjecture for affine curves. If we put together their observation and \autoref{eq}, we obtain a new proof of the following classical result, see \cite[Proposition 103]{sti13} and \cite[Proposition 250]{sti13}. 

\begin{theorem}\label{propen}
	The section map is injective for affine hyperbolic curves. The section conjecture for proper hyperbolic curves implies the section conjecture for affine hyperbolic curves.\qed
\end{theorem}

Let us show how the ideas of Borne and Emsalem fit nicely in our formalism, giving a clear picture of packets of tangential points and of the section conjecture for affine curves.

Let $U,X,D$ be as above, and define $U_{n}$ as the root stack supported over $X$ with ramification of degree $n$ along the divisor $D$, see \cite[Appendix B.2]{agv08}, \cite[§ 2.2]{be14} for the definition of root stack. The stack $U_{n}$ is an orbicurve with coarse moduli space $U_{n}\to X$ such that $U\s X$ is the schematic locus of $U_{n}$ and $U_{n}\to X$ has ramification index equal to $n$ at each point over $D$. Define
\[\hat{U}=\projlim_{n}U_{n}.\]
as the projective limit: it is an fpqc stack with natural morphisms $U\hookrightarrow\hat{U}$ and $\hat{U}\twoheadrightarrow X$. The pro-algebraic stack $\hat{U}$, called the \emph{infinite root stack}, can also be constructed using logarithmic geometry, see \cite{tv18} for details.

The natural morphism
\[\Pi_{U/k}\to\Pi_{\hat{U}/k}\simeq\projlim_{n}\Pi_{U_{n}/k}\]
is an isomorphism, this is proved in \cite[Proposition 3.2.2]{bor09}. In view of this fact, from our point of view one could simply consider $\hat{U}$ as a "complete substitute" of $U$ and decide that the section conjecture for $U$ is the section conjecture for $\hat{U}$. Let us show that this coincides with the classical approach using packets.

Fix a rational point $x\in D(k)$. Using Abhyankar's lemma, the étale fundamental gerbe $\Pi_{U^{x}/k}$ of the scheme of nearby points $U^{x}$ can be easily computed to be abelian and banded by $\hz(1)$. Let $\bar{x}$ be the closed point of $\spec(\O^{\rm h}_{X,x})$, if we consider the infinite root stack $\hat{U}^{x}$ of $\spec(\O^{\rm h}_{X,x})$ at $\bar{x}$, it is immediate to check that the structure map $\hat{U}^{x}\to\Pi_{\hat{U}^{x}/k}=\Pi_{U^{x}/k}$ induces an isomorphism between the residue gerbe $\hat{U}^{x}\times_{\spec(\O^{\rm h}_{X,x})}\bar{x}=\hat{U}^{x}_{\bar{x}}$ and $\Pi_{U^{x}/k}$. This gives us a map
\[\Pi_{U^{x}/k}=\hat{U}^{x}_{\bar{x}}\to\hat{U}_{x}=\hat{U}\times_{X}x\]
which is easily checked to be an isomorphism. It follows that the packet at $x$ identifies naturally with the isomorphism classes of rational points of $\hat{U}$ over $x$.

\begin{proof}[Proof of \autoref{propen}]
If $U$ is hyperbolic, $\rchi(U_{n})<0$ for $n$ big enough. Hence
\[U_{n}(k)\to\Pi_{U_{n}/k}(k)\]
is fully faithful for $n$ big enough, and passing to the limit the same is true for $\hat{U}$. If the section conjecture holds for proper hyperbolic curves, the section map is an equivalence for $U_{n}$ thanks to \autoref{eq}. Passing to the limit, the section map is an equivalence for $\hat{U}$, and thus the section conjecture holds for $U$.
\end{proof}

\sectionmark{The section conj. implies the hom conj.}
\section{The section conjecture implies the hom conjecture}\label{secthom}
\sectionmark{The section conj. implies the hom conj.}

If $X$ is printable, we expect the functor
\[X(T)\to\Pi_{X/k}(T)\]
to be an equivalence for a much larger class than finitely generated extensions of $k$. At least, we should have smooth schemes: we actually show that normality together with a broad finiteness condition is enough.

Recall that, for $X$ an hyperbolic curve and $T$ smooth, Mochizuki proved in \cite[Theorem A]{moc99} that $X(T)\to\Pi_{X/k}(T)$ induces a bijection between dominant morphisms $T\to X$ and sections $T\to\Pi_{X/k}$ inducing an \emph{open} homomorphism of étale fundamental groups. Here, we are concerned with \emph{all} sections $T\to\Pi_{X/k}$: the hom conjecture implies that, if $T$ is geometrically connected, a section $T\to\Pi_{X/k}$ either induces an open homomorphism of étale fundamental groups or it factorizes through $\spec k$, but no proof of this is known.

\begin{theorem}\label{point}
	Let $X$ be a smooth, proper, geometrically connected DM stack and $T$ a locally noetherian, normal scheme over $k$. Assume that, for every $t\in T$, the residue field $k(t)$ is finitely generated over $k$. If $X$ is fff, then $X(T)\to\Pi_{X}(T)$ is fully faithful. If $X$ is printable, then $X(T)\to\Pi_{X}(T)$ is an equivalence of categories.
	\begin{proof}
		Since a locally noetherian, normal scheme is a disjoint union of integral normal schemes by \cite[\href{https://stacks.math.columbia.edu/tag/033N}{Lemma 033N}]{stacks-project}, we may assume that $T$ is integral.
				
		{\bf Full faithfulness.} Let $t_{1},t_{2}:T\to X$ be two morphisms, $\pi(t_{1}),\pi(t_{2})$ their images in $\Pi_{X/k}(T)$ and $(t_{1},t_{2})\in X\times X(T)$. We have that $\uisom_{X}(t_{1},t_{2})$ is proper, unramified and hence finite over $T$ since $X$ is separated and DM. Moreover, $\Pi_{X/k}$ is a projective limit of DM, separated stacks, hence for the same reason $\uisom_{\Pi_{X/k}}(\pi(t_{1}),\pi(t_{2}))$ is profinite over $T$.
		
		Since $\uisom_{X}(t_{1},t_{2})$, $\uisom_{\Pi_{X/k}}(\pi(t_{1}),\pi(t_{2}))$ are profinite over $T$ and $T$ is integral and normal, we have that
		\[\uisom_{X}(t_{1},t_{2})(T)=\uisom_{X}(t_{1},t_{2})(k(T)),\]
		\[\uisom_{\Pi_{X/k}}(\pi(t_{1}),\pi(t_{2}))(T)=\uisom_{\Pi_{X/k}}(\pi(t_{1}),\pi(t_{2}))(k(T)),\]
		and hence
		\[\uisom_{X}(t_{1},t_{2})(T)\xar{\sim}\uisom_{\Pi_{X/k}}(\pi(t_{1}),\pi(t_{2}))(T)\]
		since by hypothesis
		\[\uisom_{X}(t_{1},t_{2})(k(T))\xar{\sim}\uisom_{\Pi_{X/k}}(\pi(t_{1}),\pi(t_{2}))(k(T)).\]
		
		{\bf Essential surjectivity.} Fix a morphism $\tau:T\to\Pi_{X/k}$. Consider the "universal cover" $\tilde{X}=X\times_{\Pi_{X/k}}T\to X\times T$, which is an algebraic space thanks to \autoref{cover}: our hypothesis that $X$ is printable implies that, for every finitely generated extension $k'/k$, $\tilde{X}(k')\to T(k')$ is bijective.
		
		Let $R$ be a DVR over $k$ with fraction and residue fields finitely generated over $k$, and suppose we have a morphism $\spec R\to T$. Since $k(R)$ is finitely generated over $k$, $\spec k(R)\to T$ lifts uniquely to $\tilde{X}$. The morphism $\tilde{X}\to T$ is separated and universally closed, since $\tilde{X}\to X\times T$ is profinite (it can be obtained by base change from the diagonal of $\Pi_{X/k}$) and $X\times T\to T$ is proper. Thanks to \cite[\href{https://stacks.math.columbia.edu/tag/0A3X}{Lemma 0A3X, Lemma 0A3W, Lemma 03KU}]{stacks-project} this is enough to apply the valuative criterion, i.e. we get that $\spec R\to T$ lifts uniquely to $\tilde{X}$. We thus obtain that $\tilde{X}(R)\to T(R)$ is bijective, too. If $A$ is either a field finitely generated over $k$ or a DVR with fraction and residue fields finitely generated over $k$, call 
		\[\iota:T(A)\to\tilde{X}(A)\]
		the inverse map.
				
		Thanks to \autoref{cover}, there exists a finite étale gerbe $\Phi$ and a representable morphism $X\to\Phi$, the fiber product $X\times_{\Phi}T$ is an algebraic space. The natural morphism $\Pi_{X/k}\to\Phi$ induces a morphism 
		\[\omega:\tilde{X}\to X\times_{\Phi}T.\]
		The reader may keep in mind the particular case in which $X$ is an algebraic space and $\Phi=\spec k$: we use the morphism $X\to\Phi$ only to "kill" the inertia of $X\times T$.
		
		Consider now the generic point $\xi:\spec k(T)\to T$, and let $S\s X\times_{\Phi}T$ be the closure of $\omega(\iota(\xi)):\spec k(T)\to X\times_{\Phi}T$ with the reduced structure, and $p:S\to T$ the projection. We have that $S$ is an integral algebraic space. The situation is illustrated in the following diagram.
		
		\[\begin{tikzcd}[row sep=large]
			\tilde{X}\ar[rr, bend left,"\omega"]	&	S\rar[hook]\ar[d,"p"]		&	X\times_{\Phi}T\dar		\\
			\spec k(T)\ar[r,"\xi"]\uar[bend left,pos=0.45,"\iota(\xi)"]\ar[ur]		&	T\ar[r,equal]\ar[u,dotted,bend left,"\exists ?"]			&	T
		\end{tikzcd}\]
		
		Observe that it is sufficient to prove that $p:S\to T$ is an isomorphism. In fact, if $S\to T$ is an isomorphism, we have an induced morphism $x:T\to X\times_{\Phi} T\to X$ such that $\pi(x):T\to\Pi_{X/k}$ is generically isomorphic to the original morphism $\tau:T\to \Pi_{X/k}$. Thanks to what we have shown in the preceding point, the fact that $\tau$ and $\pi(x)$ are generically isomorphic implies that they are isomorphic, hence $\tau$ is in the essential image of $X(T)\to\Pi_{X/k}(T)$. Let us show that $p:S\to T$ is an isomorphism.
		
		{\bf Step 1:} the map $S(k')\to T(k')$ is injective for any field extension $k'/k$. Since $S$ is of finite type over $k$ and all the points of $T$ have residue field finitely generated over $k$, it is enough to do so for $k'/k$ finitely generated. For every point $s:\spec k'\to S$ with $k'/k$ finitely generated, consider $\omega(\iota(p(s))):\spec k'\to X\times_{\Phi}T$. It is enough to show that 
		\[\omega(\iota(p(s)))=s,\]
		and it is enough to do so for $k'=k(s)$. We prove this by induction on the Krull height of $s$ in $S$ (the height can be defined by passing to an étale neighbourhood which is a scheme).
		
		If $s$ has height $0$, then $s$ is the generic point $\omega(\iota(\xi))$ of $S$ and $p(s)$ maps to the generic point $\xi$ of $T$, thus $k(T)= k(s)$ and $s=\omega(\iota(\xi))=\omega(\iota(p(s)))$.
		
		If $\on{ht}_{S}(s)> 0$, there exists a germ of a non-constant curve on $S$ passing through $s$. More precisely, there exists a noetherian DVR $R$ with fraction and residue fields finitely generated over $k$ and a morphism $r:\spec R\to S$ such that the closed point maps to $s$ and the open point maps to a point $s_{0}\neq s$. In order to find $R$, take an étale neighbourhood $(S',s')$ of $s$ which is a scheme, and choose $R$ as the normalization of a dimension $1$ integral quotient of $\O_{S',s'}$.
		
		Now consider $\omega(\iota(r)):\spec R\to X\times_{\Phi}T$. We have $\omega(\iota(p(r)))_{k(R)}=r_{k(R)}$ by induction hypothesis, this implies that $\omega(\iota(p(r)))=r$ since $S\to T$ is separated. Let $R/\mf{m}$ be the residue field of $R$, we have $k(s)\s R/\mf{m}$, 
		\[\omega(\iota(p(s)))_{R/\mf{m}}=\omega(\iota(p(r)))_{R/\mf{m}}=r_{R/\mf{m}}=s_{R/\mf{m}}\]
		and thus $\omega(\iota(p(s)))=s$. This concludes step 1.
		
		{\bf Step 2:} the map $S(k')\to T(k')$ is bijective. We already know that it is injective, moreover $T$ is integral and $S\to T$ is proper, thus $S\to T$ is bijective set-theoretically. Let $t\in T$ be any point, we thus know that $S_{t}$ has exactly one point $s\in S_{t}$. Since $S_{t}$ is of finite type over $k(t)$, it follows that $k(s)$ is finite over $k(t)$. We are in characteristic $0$ and thus $k(s)$ is separable over $k(t)$, since $S_{t}(k')\to \spec k(t)(k')$ is injective for every extension $k'/k$ we get that $k(s)=k(t)$, and this concludes step 2.
		
		{\bf Step 3:} $p:S\to T$ is an isomorphism. Thanks to the previous steps, $S$ is quasi-finite over $T$ and thus it is a scheme, see \cite[\href{https://stacks.math.columbia.edu/tag/03XX}{Proposition 03XX}]{stacks-project}.
		
		If $T$ is quasi-compact and quasi-separated, we can apply Zariski's main theorem \cite[\href{https://stacks.math.columbia.edu/tag/05K0}{Lemma 05K0}]{stacks-project} and there exists a factorization $S\to T'\to T$ with $S\to T'$ an open immersion and $T'\to T$ finite. Since $S$ is integral, we may assume $T'$ integral too and $k(S)=k(T')=k(T)$. Since $T$ is normal, $T'\to T$ is an isomorphism. It follows that $S\to T$ is a bijective open immersion, i.e. an isomorphism.
		
		If $T$ is not quasi-compact and quasi-separated, cover it by open affine schemes $T_{i}$ with restrictions $S_{i}\to T_{i}$. For each $i$, the argument above works since $T_{i}$ is quasi-compact and quasi-separated, hence $S_{i}\to T_{i}$ is an isomorphism. It follows that $S\to T$ is an isomorphism.
	\end{proof}
\end{theorem}

\begin{remark}
	We cannot hope to remove completely the normality hypothesis from \autoref{point}. Consider an integral, projective curve $X$ of geometric genus at least $2$ with only a cuspidal singularity and smooth normalization $\overline{X}$ over a number field $k$. Passing to $\C$, we may check that $\Pi_{X/k}=\Pi_{\overline{X}/k}$. If the section conjecture holds, $\overline{X}$ is printable: \autoref{point} without the normality hypothesis would give us a section $X\to\overline{X}$.
\end{remark}

\begin{corollary}\label{curvesimpl}
	If smooth, proper, hyperbolic curves satisfy the section conjecture, then they satisfy the hom conjecture.
	\begin{proof}
		If hyperbolic curves satisfy the section conjecture, then they are printable thanks to \autoref{curvessect}. Hence, they satisfy the hom conjecture thanks to \autoref{point}.
	\end{proof}
\end{corollary}

Thanks to \autoref{curvesimpl}, we can also see the anabelian conjecture proved by Mochizuki as a particular case of the section conjecture, rather than a different one.

\autoref{point} allows us to prove that the topological fundamental group of a printable DM stack has no abelian finite index subgroup. We know no other result of the form "if a variety shows anabelian behaviour, then its fundamental group is far from being abelian": conjectures and theorems are always in the other direction.

\begin{proposition}\label{anvsab}
	Let $X$ be a printable DM stack of positive dimension. Then $\pi_{1}(X_{\bar{k}})$ has no finite index abelian subgroups.
	\begin{proof}
		Thanks to \autoref{cover} and \autoref{etale}, up to a finite extension of $k$ and a finite étale covering of $X$ we may assume that $X$ is an algebraic space. 
		
		Assume by contradiction that a finite index abelian subgroup exists, up to another a finite extension of $k$ and finite étale covering of $X$ we may assume that $\pi_{1}(X_{\bar{k}})$ is abelian and $X$ has a rational point $x_{0}\in X(k)$. Let $\rm{Sm}_{k}$ be the category of smooth varieties over $k$. Since $X$ is printable, thanks to \autoref{point} $X$ and $\Pi_{X/k}$ define two naturally equivalent functors $\rm{Sm}_{k}^{\op}\to\set$ (by taking equivalence classes of $\Pi_{X/k}(T)$ for every $T\in\rm{Sm}_{k}$). The fact that the fundamental group of $X_{\bar{k}}$ is abelian implies that the gerbe $\Pi_{X/k}$ is abelian and hence its functor is enriched in groups with identity $\pi(x_{0})\in\Pi_{X/k}(x_{0})$, thus the same is true for the functor defined by $X$ and $x_{0}$.
		
		Now take an étale cover $U\to X$ with $U$ a scheme, and let $R=U\times_{X}U$. Since $U$ and $R$ are smooth varieties, $X(U)$ and $X(R)$ are groups with the structure inherited from $\Pi_{X/k}(U)$ and $\Pi_{X/k}(R)$, this allows us to construct the usual maps $m:X\times X\to X$, $i:X\to X$ giving the group structure to $X$. Hence, the functor of points of $X$ is enriched in groups over the whole category of schemes over $k$ and not just the smooth ones. This implies that $X$ is not only an algebraic space but also a scheme, see \cite[Theorem 4.1]{art69}.
		
		Hence, $X$ is actually a proper group scheme, i.e. an abelian variety. But it is well known that an abelian variety of positive dimension is not printable, see for instance \href{https://mathoverflow.net/questions/92927/}{MathOverflow 92927} where a proof is given for elliptic curves (the proof actually works without modifications for positive dimensional abelian varieties).
	\end{proof}
\end{proposition}

\section{Elementary anabelian stacks}\label{elemsect}

Recall that Grothendieck defined a geometrically connected variety $X$ as elementary anabelian if there exists a chain of morphisms
\[X=X_{N}\to X_{N-1}\to\dots X_{1}\to X_{0}=\spec k\]
such that $X_{i+1}\to X_{i}$ is a smooth fibration whose fibers are hyperbolic curves, see \cite{gro97}. We want to extend this definition to \emph{elementary anabelian stacks}.

\begin{definition}
	Let $X,Y$ be DM stacks. A morphism $Y\to X$ is a \emph{family of orbicurves} if it is smooth, proper, and its fibers are geometrically connected orbicurves.
\end{definition}

\begin{definition}\label{elandef}
	Let $k$ be a field of characteristic $0$. A DM stack is \emph{constructible by fibrations} over $k$ if it can be constructed by recursion in the following way.
	\begin{enumerate}
		\item $\spec k$ is constructible by fibrations.
		\item If $Y\to X$ is a family of hyperbolic orbicurves and $X$ is constructible by fibrations, then $Y$ is constructible by	fibrations.
		\item If $Y\to X$ is finite, representable and étale, then $X$ is constructible by fibrations if and only if $Y$ is constructible by fibrations.
	\end{enumerate}
	We say that $X$ is an \emph{elementary anabelian stack} if there exists a field extension $k'/k$ such that $X_{k'}$ is constructible by fibrations.
\end{definition}

\begin{lemma}\label{eabc}
	Let $X$ be a DM stack over a field $k$ of characteristic $0$, and $K/k$ any extension. Then $X$ is elementary anabelian if and only if $X_{K}$ is elementary anabelian.
	\begin{proof}
		If $X$ is elementary anabelian, there exists an extension $k'/k$ such that $X_{k'}$ is constructible by fibrations. It is possible to find an extension $K'/k$ containing both $k'$ and $K$ as subextensions. It follows that $X_{K'}$ is constructible by fibrations and thus $X_{K}$ is elementary anabelian. The other implications is trivial.
	\end{proof}
\end{lemma}

\begin{lemma}\label{cfvsea}
	Let $X$ be an elementary anabelian stack over a field $k$ of characteristic $0$. There exists a diagram of field extensions
	\[\begin{tikzcd}[row sep=tiny]
							&	k'	&								&	\C	\\
			k\ar[ur,hook]	&		&	h\ar[ul,hook]\ar[ur,hook]
	\end{tikzcd}\]
	such that the extensions $k'/k$, $h/\Q$ are finitely generated, and $X_{k'}$ descends to a DM stack $Y$ over $h$ which is constructible by fibrations.

	If $k$ is finitely generated over $\Q$ and $k\s\C$, we may choose $k'=h$, i.e. there exists a finitely generated subextension $\C/k'/k$ such that $X_{k'}$ is constructible by fibrations.
	\begin{proof}
		The first part follows from standard arguments about finite presentation. For the second part, find $k'$ as before and embed it in $\C$ as an extension of $k$ using the fact that $\C$ is algebraically closed of infinite transcendence degree over $k$.
	\end{proof}
\end{lemma}

\begin{corollary}\label{eac}
	A DM stack over $\C$ is an elementary anabelian stack if and only if it is constructible by fibrations.
	\begin{proof}
		The "if" part is by definition. Let $X$ be elementary anabelian over $\C$, it descends to a DM stack $Y$ over some finitely generated subfield $k\s \C$. Apply \autoref{cfvsea} to $Y$, we find a finitely generated subextension $\C/k'/k$ such that $Y_{k'}$ is constructible by fibrations. It follows that $Y_{\C}=X$ is constructible by fibrations.
	\end{proof}
\end{corollary}

\subsection{Topology of elementary anabelian stacks}

It will be useful to define and study the topological counterparts of elementary anabelian stacks. Mainly, we do so because we will need the long exact sequence of a fibration: in the algebraic setting, the standard reference \cite{fri73} does not cover stacks. As a workaround, we will pass through topology and then get back to the algebraic setting using Serre's good groups. We refer to \cite{noo05} and \cite{noo14} for the theory of topological stacks and orbifolds, fibrations and long exact homotopy sequences.

\begin{definition}
	A \emph{topological orbicurve} is a complete orbifold $X$ of dimension $2$ with a coarse moduli space $X\to S$ which restricts to an isomorphism on the complementary of a finite subset of $S$ and such that $S$ is a compact, orientable surface.
\end{definition}

Topological orbicurves have a rational Euler characteristic analogously to algebraic orbicurves, thus we may define hyperbolic, elliptic and parabolic topological orbicurves.

\begin{definition}
	Elementary anabelian topological stacks are topological DM stacks defined by recursion in the following way.
	\begin{enumerate}
		\item The point is elementary anabelian.
		\item If $Y\to X$ is a fibration whose fibers are hyperbolic topological orbicurves and $X$ is elementary anabelian, then $Y$ is elementary anabelian.
		\item If $Y\to X$ is a finite covering space, then $X$ is elementary anabelian if and only if $Y$ is elementary anabelian.
	\end{enumerate}
\end{definition}

In order to pass from the topological to the algebraic setting, we need to check that all groups involved are good in the sense of Serre. Recall that a discrete group $G$ is good in the sense of Serre if the natural homomorphism
\[\H^{q}(\hat{G},M)\to\H^{q}(G,M)\]
is an isomorphism for every finite $G$-module $M$, where $\hat{G}$ is the profinite completion of $G$. The reason why we are interested in good groups is the following: if $G/K=H$ is an extension with $G$ good and $H$ finitely generated, then
\[1\to \hat{K}\to \hat{G}\to\hat{H}\to 1\]
is exact, see \cite[§ I.2.6 Exercises 1,2]{ser94}.

Finite groups are obviously good, and fundamental groups of compact, orientable surfaces are known to be good: this was already known to Serre when he introduced the concept, see \cite[p. 138]{scm04}. In \cite[Proposition 3.7]{gjz08}, it is proved that fundamental groups of topological orbicurves are good.

We need now to strengthen the concept of good groups in order to obtain a definition which is stable under extension. We will do so by requesting that our groups are of type $\FLi$. Recall that a group $G$ is of type $\FLi$ if $\Z$ as a $\Z[G]$-module has a resolution by free $\Z[G]$-modules of finite rank, see \cite[Chapter VIII]{bro94}. Groups of type $\FLi$ are finitely generated, see \cite[§ VIII.4 Exercise 1]{bro94}. If a group $G$ is of type $\FLi$ and $M$ is a finite $G$-module, then $\H^{q}(G,M)=\on{Ext}^{q}_{\Z[G]}(\Z,M)$ is finite since it is torsion and we can use the finite rank resolution of $\Z$ to show that it is finitely generated.

Finite groups are of type $\FLi$ since the resolution can easily be constructed step by step. Moreover, if there exists a CW complex $X$ of type $K(G,1)$ whose $n$-skeleton is finite for every $n$, the cellular chain complex of the universal covering of $X$ shows that $G$ is of type $\FLi$, see \cite[§ I.4]{bro94}. In particular, the fundamental group of surfaces is of type $\FLi$. 

\begin{definition}
	We say that a group is \emph{very good} if it is good and of type $\FLi$.
\end{definition}

\begin{lemma}\label{vgext}
	Let $H=G/K$ be an extension of groups. If $H,K$ are very good, then $G$ is very good.
	\begin{proof}
		If $H,K$ are of type $\FLi$, it is known that $G$ is of type $\FLi$, see \cite{wal61}.
		
		If $H,K$ are good, a sufficient condition for $G$ to be good is that $K$ is finitely generated and $\H^q(K,M)$ is finite for every finite $G$-module $M$, see \cite[§ I.2.6 Exercises 1,2]{ser94}. If $K$ is of type $\FLi$, this is automatically satisfied.
	\end{proof}  
\end{lemma}

Recall that two groups are called commensurable if they have isomorphic finite index subgroups.

\begin{lemma}\label{vgcomm}
	If $H,G$ are commensurable and $G$ is very good, then $H$ is very good.
	\begin{proof}
		The fact that $H$ is good is proved in \cite[Lemma 3.2]{gjz08}. We want to show that it is $\FLi$, too.
		
		Suppose first that $H$ is a finite index subgroup of $G$. Since $\Z[G]$ is a free $\Z[H]$ module of finite rank, it is immediate to check that $H$ is of type $\FLi$.
		
		Suppose now that $G$ is a finite index subgroup of $H$. By the preceding case, we may reduced to the case in which $G$ is normal of finite index in $H$. The thesis then follows from \autoref{vgext}.
		
		The general case directly follows from the first two.
	\end{proof}
\end{lemma}

\begin{lemma}\label{elemtop}
	Let $X$ be an elementary anabelian topological stack. Then
	\begin{enumerate}
		\item $X$ is of type $K(G,1)$,
		\item $X$ is uniformizable, i.e. there exists a finite covering space $Y\to X$ with $Y$ a manifold,
		\item $\pi_{1}(X)$ is very good and residually finite.
	\end{enumerate}
	\begin{proof}
		Let us prove this by recursion using the rules defining elementary anabelian topological stacks.
		
		If $X$ is just a point, then it clearly satisfies the thesis.
		
		Let $Y\to X$ is a finite covering space. Clearly, $Y$ is of type $K(G,1)$ and uniformizable if and only if the same holds for $X$. By \autoref{vgcomm}, $\pi_{1}(X)$ is very good if an only if the same holds for $\pi_{1}(Y)$. If $H=G/K$ is a finite group, then it is immediate to show that $G$ is residually finite if and only if $K$ is residually finite. Using a Galois closure $Z\to Y\to X$, we get that $\pi_{1}(X)$ is residually finite if and only if the same holds for $\pi_{1}(Y)$.
		
		If $X$ is a complete, orientable surface of genus at least $2$, it is clearly of type $K(G,1)$ and uniformizable. The fact that $\pi_{1}(X)$ is good is proved in \cite[Proposition 3.7]{gjz08}. Since $X$ is of type $K(G,1)$ and we can endow it with a structure of finite CW complex, it follows that $\pi_{1}(X)$ is of type $\FLi$, and thus very good. The fact that $\pi_{1}(X)$ is residually finite is well-known, see \cite{hem72}.
		
		Let $X$ be an hyperbolic topological orbicurve. Thanks to \cite[Proposition 5.1, Corollary 5.9]{bn06}, there exists some finite covering space $Y\to X$ with $Y$ an hyperbolic surface, thus the thesis follows from the preceding case.
		
		Now suppose that $X$ satisfies the thesis and let $Y\to X$ be a fibration whose fibers are hyperbolic topological orbicurves, and let $F$ be a fiber. We have a commutative diagram
		\[\begin{tikzcd}
			1\rar	&	\pi_{1}(F)\rar\dar		&	\pi_{1}(Y)\rar\dar		&	\pi_{1}(X)\rar\dar		&	1	\\
			1\rar	&	\hat{\pi_{1}(F)}\rar	&	\hat{\pi_{1}(Y)}\rar	&	\hat{\pi_{1}(X)}\rar	&	1
		\end{tikzcd}\]
		where the first row is exact being the long exact sequence of the fibration. It follows that $Y$ is of type $K(G,1)$. Since $\pi_{1}(X)$ and $\pi_{1}(F)$ are very good, $\pi_{1}(Y)$ is very good too thanks to \autoref{vgext} and the second row is exact. Moreover, the left and right vertical arrows are injective since $\pi_{1}(F),\pi_{1}(X)$ are residually finite, thus the central vertical arrow is injective by diagram chasing and hence $\pi_{1}(Y)$ is residually finite.
		
		It remains to show that $Y$ is uniformizable. Since $X$ is uniformizable, up to passing to a finite covering space we may assume that $X$ is a manifold. We know that $F$ is uniformizable, let $F'\to F$ be a finite covering space which is a surface, and let $Q\s\hat{\pi_{1}(F)}$ be the set theoretic complement of $\hat{\pi_{1}(F')}$. Since $\pi_{1}(F')$ has finite index in $\pi_{1}(F)$, we have that $Q\s\hat{\pi_{1}(F)}$ is compact. By injectivity of $\hat{\pi_{1}(F)}\to\hat{\pi_{1}(Y)}$, the image of $Q$ in $\hat{\pi_{1}(Y)}$ does not contain $1$. By compactness of $Q$ it follows that there exists a finite quotient $\hat{\pi_{1}(Y)}\to G$ such that the image of $Q$ in $G$ does not contain $1$. Let $Y'\to Y$ be the associated $G$-covering. I claim that $Y'$ is a manifold.
		
		In fact, let $K\s\hat{\pi_{1}(F)}$ be the kernel of the composition $\hat{\pi_{1}(F)}\to\hat{\pi_{1}(Y)}\to G$. Since the image of $Q$ in $G$ does not contain $1$, then $K\s \hat{\pi_{1}(F')}$, hence the associated (possibly disconnected) $G$-covering $F''\to F$ dominates $F'$ and is thus a manifold. In particular, $Y'\to X$ is a fibration whose base and fibers are manifolds, and thus it is a manifold too. 
	\end{proof}
\end{lemma}

\begin{corollary}
	The étale homotopy type of an elementary anabelian stack over an algebraically closed field of characteristic $0$ is of type $K(G,1)$.
	\begin{proof}
		Thanks to \autoref{cfvsea}, we may reduce to the case in which $k=\C$ and $X$ is constructible by fibrations. Let $X$ be an elementary anabelian stack, the associated topological stack $X^{\rm an}$ in the sense of \cite{noo05} is an elementary anabelian topological stack. Now apply \cite[Theorem 6.7]{am69} and \autoref{elemtop}.
	\end{proof}
\end{corollary} 

\subsection{An additional property of elementary anabelian stacks}

The class of elementary anabelian stacks is by definition stable under a certain number of operations. There is another natural, non-obvious operation under which they are stable: if $Y\to X$ is an elementary anabelian stack and $Y\to X$ is a proper étale morphism, for instance a proper étale gerbe, then $Y$ is elementary anabelian.

We remark that this property is not strictly necessary for the main purposes of the paper: it just seems appropriate to highlight it while introducing elementary anabelian stacks.

We refer to \cite[\href{https://stacks.math.columbia.edu/tag/06QB}{Section 06QB}]{stacks-project} for the definition and main properties of gerbes.

\begin{lemma}\label{peger}
	Let $Y\to X$ be a morphism of finite type of algebraic stacks which is a gerbe. The following are equivalent.
	\begin{enumerate}
		\item $Y\to X$ is proper étale.
		\item $Y\to X$ is separated and DM.
		\item The diagonal $Y\to Y\times_{X}Y$ is finite étale. 
		\item The relative inertia $I_{Y/X}\to Y$ is finite étale. 
	\end{enumerate}
	\begin{proof}
		(1) $\Rightarrow$ (2). By definition.
		
		(2) $\Rightarrow$ (3). Since $Y\to X$ is separated and DM, the diagonal is proper and unramified. Since the diagonal of a gerbe is always flat \cite[\href{https://stacks.math.columbia.edu/tag/0CPR}{Lemma 0CPR}]{stacks-project}, we conclude.

		(3) $\Rightarrow$ (1).	Since $Y\to X$ is gerbe, it is a universal homeomorphism by \cite[\href{https://stacks.math.columbia.edu/tag/06R9}{Lemma 06R9}]{stacks-project}, it is separated since the diagonal is finite and it is of finite type by hypothesis, thus it is proper. Moreover, $Y\to X$ is flat since it is a gerbe \cite[\href{https://stacks.math.columbia.edu/tag/0CPS}{Proposition 0CPS}]{stacks-project} and unramified since it is of finite type with étale diagonal \cite[\href{https://stacks.math.columbia.edu/tag/0CJ0}{Lemma 0CJ0}]{stacks-project}, thus we conclude that it is étale, too.
		
		(3) $\Leftrightarrow$ (4). Since $Y\to X$ is gerbe, the diagonal is an fppf covering by \cite[\href{https://stacks.math.columbia.edu/tag/0CPS}{Proposition 0CPS}]{stacks-project}. It follows that the diagonal is finite étale if and only if the same holds for the relative inertia.
	\end{proof}
\end{lemma}

If $\Phi\to\spec k$ is a proper étale gerbe over a field $k$, then it is immediate to check that it is a finite étale gerbe in the sense of \cite{bv15}. 

\begin{lemma}\label{stein}
	Let $f:Y\to X$ be a proper étale morphism of DM stacks locally of finite type with Stein factorization
	\[Y\to\spec f_{*}\O_{Y}\to X.\]
	Then $Y\to \spec f_{*}\O_{Y}$ is a proper étale gerbe and $\spec f_{*}\O_{Y}\to X$ is representable finite étale.
	\begin{proof}
		Up to passing to an étale cover and since $f_{*}$ commutes with étale base change, we may assume that $X$ is a scheme. If $X$ is a scheme, the inertia stack $I_{Y}$ of $Y$ coincides with the relative inertia stack $I_{Y/X}$, see \cite[\href{https://stacks.math.columbia.edu/tag/04Z6}{Lemma 04Z6}]{stacks-project}. Since $Y\to X$ is étale, it follows that $I_{Y}\simeq I_{Y/X}$ is étale over $Y$ thanks to \cite[\href{https://stacks.math.columbia.edu/tag/0CJ0}{Lemma 0CJ0}]{stacks-project}. In particular, $I_{Y}$ is flat over $Y$, and thus $Y$ is a gerbe over its coarse moduli sheaf $M$ which is an algebraic space too, see \cite[\href{https://stacks.math.columbia.edu/tag/06QJ}{Proposition 06QJ}]{stacks-project} and \cite[\href{https://stacks.math.columbia.edu/tag/06QD}{Lemma 06QD}]{stacks-project}.
		
		Since $f$ is proper and $X$ is locally of finite type, pushforward of coherent sheaves is coherent, see \cite{fal03}, and hence $\spec f_{*}\O_{Y}\to X$ is a finite morphism. We have a natural morphism $M\to \spec f_{*}\O_{Y}$ since $X$, and thus $\spec f_{*}\O_{Y}$, is a scheme. On the other hand, $M\to X$ is proper and quasi-finite, hence affine, and this gives us a natural morphism in the other direction $\spec f_{*}\O_{Y}\to M$. These are easily checked to be inverses.
		
		Thus, we know that $Y\to M=\spec f_{*}\O_{Y}$ is a gerbe and $M\to X$ is finite. Since $Y\to X$ is separated and DM and $M\to X$ is representable finite, by \cite[\href{https://stacks.math.columbia.edu/tag/050M}{Lemma 050M}]{stacks-project} we get that $Y\to M$ is separated and DM and thus proper étale by \autoref{peger}.	It remains to prove that $M\to X$ is étale, and this follows from the fact that both $Y\to M$, $Y\to X$ are surjective étale.
	\end{proof}
\end{lemma}

\begin{lemma}\label{elemgerbe}
	If $X$ is an elementary anabelian stack and $f:Y\to X$ is a proper étale morphism, then $Y$ is an elementary anabelian stack.
	\begin{proof}
		Thanks to \autoref{cfvsea}, we may reduce to the case in which $k=\C$ and $X$ is constructible by fibrations, and we want to show that $Y$ is constructible by fibrations too. Thanks to \autoref{stein}, we may furthermore reduce to the case in which $Y\to X$ is a proper étale gerbe.
		
		Consider a geometric point $y\in Y(\C)$ and its image $x\in X(\C)$. The fiber $Y_{x}$ is a gerbe of the form $BG$ for some finite group $G$. Since $Y\to X$ is a proper étale gerbe, the diagonal $Y\to Y\times_{X}Y$ is finite étale. Passing to the associated topological stacks in the sense of \cite{noo05}, \cite{noo14}, this tells us that $Y^{\rm an}\to X^{\rm an}$ is a fibration with fiber $BG^{\rm an}$, and we may thus consider the topological homotopy exact sequence
		\[1\to G\to \pi_{1}^{\rm top}(Y^{\rm an})\to\pi_{1}^{\rm top}(X^{\rm an})\to 1,\]
		where $\pi_{2}^{\rm top}(X)$ is $0$ by \autoref{elemtop}. Since $G$, $\pi_{1}^{\rm top}(X)$ are very good, we can pass to profinite completions 
		\[1\to G\to\hat{\pi_{1}^{\rm top}(Y^{\rm an})}=\pi_{1}(Y)\to\hat{\pi_{1}^{\rm top}(X^{\rm an})}=\pi_{1}(X)\to 1.\]
		Since $G$ is finite, there exists a connected, finite étale cover $Z\to Y$ such that $\pi_{1}(Z)\cap G=\{1\}\s\pi_{1}(Y)$.
		
		Consider now the composition $g:Z\to Y\to X$, it is a proper étale morphism. I claim that it is representable. Consider the Stein factorization
		\[Z\to\spec g_{*}\O_{Z}\to X,\]
		it is enough to show that $Z\to \spec g_{*}\O_{Z}$ is an isomorphism. By \autoref{stein}, we know that $Z\to \spec g_{*}\O_{Z}$ is a proper étale gerbe and that $\spec g_{*}\O_{Z}\to X$ is a finite étale cover. Take any geometric fiber of $Z\to\spec g_{*}\O_{Z}$, it has the form $BH$ for some finite group $H$, we want to show that $H$ is trivial.
		
		By the same argument as above, we get an embedding $H\s\pi_{1}(Z)\s\pi_{1}(Y)$, and by construction $H$ maps to the identity in $\pi_{1}(X)$. Since $\pi_{1}(Z)$ intersects the kernel of $\pi_{1}(Y)\to\pi_{1}(X)$ trivially, it follows that $H$ is trivial too.
		
		Hence, we have two finite étale covers $Z\to Y$ and $Z\to X$: since $X$ is constructible by fibrations, $Z$ and $Y$ are constructible by fibrations too.
	\end{proof}
\end{lemma}

\subsection{From curves to elementary anabelian stacks}

\begin{theorem}\label{elem}
		Elementary anabelian stacks over a field $k$ finitely generated over $\Q$ are fff. 
		
		If the section conjecture holds for smooth, proper, hyperbolic curves defined over fields finitely generated over $\Q$, then elementary anabelian stacks defined over fields finitely generated over $\Q$ are printable.
	\begin{proof}
		Thanks to \autoref{anabext} and \autoref{cfvsea}, it is enough to prove the theorem for DM stacks constructible by fibrations. We do this only for printability, the argument for fff is analogous. We are going to check that printability is preserved along the elementary operations that define DM stacks constructible by fibrations. We may assume that smooth, proper, hyperbolic orbicurves are printable thanks to \autoref{eq}.
		
		Obviously, $\spec k$ is printable since $\Pi_{\spec k/k}=\spec k$. If $Y\to X$ is finite étale, then by \autoref{etale} $Y$ is printable if and only if $X$ is printable. We only have to check that printability is preserved along families of hyperbolic orbicurves.
		
		Let $Y\to X$ be a family of hyperbolic orbicurves, and assume that $X$ is printable. Call $\Pi_{Y/X}$ the fiber product $X\times_{\Pi_{X/k}}\Pi_{Y/k}$, we have a natural 2-commutative diagram
		\[\begin{tikzcd}
			Y\ar[r]\ar[dr]	&	\Pi_{Y/X}\dar\rar	&	\Pi_{Y/k}\dar	\\
							&	X\rar				&	\Pi_{X/k}
		\end{tikzcd}\]
		For any extension $k'/k$ and morphism $x:\spec k'\to X$, consider the fiber
		\[\Pi_{Y/X,x}=\Pi_{Y/X}\times_{X}\spec k'=\Pi_{Y/k}\times_{\Pi_{X/k}}\spec k'.\]
		There is a natural map $Y_{x}\to\Pi_{Y/X,x}$.
		
		\emph{Claim:} $Y_{x}\to\Pi_{Y/X,x}$ is the étale fundamental gerbe of $Y_{x}$. Thanks to \autoref{basechange}, we may assume $k'=k=\bar{k}$ is algebraically closed. Fix a base point $y\in Y_{x}$. Then, since $X$ has trivial topological second homotopy group, there is an exact sequence of étale fundamental groups
		\[0\to\pi_{1}^{\rm top}(Y_{x},y)\to\pi_{1}^{\rm top}(Y,y)\to\pi_{1}^{\rm top}(X,x)\to 0.\]
		Since $\pi_{1}^{\rm top}(X,x)$ is good in the sense of Serre thanks to \autoref{elemtop}, we may pass to profinite completions, i.e. étale fundamental groups:
		\[0\to\pi_{1}(Y_{x},y)\to\pi_{1}(Y,y)\to\pi_{1}(X,x)\to 0.\]
		Since $\Pi_{Y/X,x}=\Pi_{Y/k}\times_{\Pi_{X/k}}\spec k(x)$, there is also a short exact sequence 
		\[0\to\aut_{\Pi_{Y/X,x}}(y)\to\aut_{\Pi_{Y/k}}(y)\to\aut_{\Pi_{X/k}}(x)\to 0,\]
		and there are natural identifications
		\[\pi_{1}(Y_{x},y)=\aut_{\Pi_{Y_{x}/k}}(y),~\pi_{1}(Y,y)=\aut_{\Pi_{Y/k}}(y),~\pi_{1}(X,x)=\aut_{\Pi_{X/k}}(x).\]
		These fit in a commutative diagram of short exact sequences, identifying $Y_{x}\to\Pi_{Y/X,x}$ with the étale fundamental gerbe $Y_{x}\to\Pi_{Y_{x}/k}$.
		
		So we know that $\Pi_{Y/X,x}$ is the étale fundamental gerbe of $Y_{x}$. Since we are assuming that hyperbolic orbicurves are printable, we get that $Y\to\Pi_{Y/X}$ is an equivalence for every finitely generated extension $k'/k$ by working fiberwise. Moreover, we are assuming that $X(k')\to\Pi_{X/k}(k')$ is an equivalence, thus the same holds for its base change $\Pi_{Y/X}(k')\to\Pi_{Y/k}(k')$. It follows that the composition $Y(k')\to \Pi_{Y/k}(k')$ is an equivalence, as desired.
	\end{proof}
\end{theorem}

Finally, let us give a version of \autoref{elem} for classical elementary anabelian varieties: for them, we can slightly relax the hypotheses since we can work with a fixed base field. Varying the base field is necessary in the passage from curves to orbicurves, but not in the one from curves to elementary anabelian varieties.

\begin{theorem}\label{classelem}
	Let $k$ be a finitely generated extension of $\Q$. If the section conjecture holds for smooth, proper, hyperbolic curves over $k$, then it holds for proper elementary anabelian varieties over $k$. Moreover, elementary anabelian varieties respect the injectivity part of the section conjecture.
	\begin{proof}
		If we assume that the section conjecture holds for an elementary anabelian variety $X/k$ and for all smooth, proper hyperbolic curves over $k$, and $Y\to X$ is a fibration in hyperbolic curves, then we can repeat the argument contained in the proof of \ref{elem} in order to show that the section conjecture holds for $Y/k$. The injectivity part is analogous.
	\end{proof}
\end{theorem}


\appendix
\section{Étale fundamental gerbes}\label{etaleapp}

Almost everything in this appendix is already known to the mathematical community, we claim no originality. In particular, most of the ideas and results are already implicit in \cite{bv15} and in the original paper by Deligne \cite{del89}. Anyway, we could not find a satisfying reference, since \cite{bv15} is mostly concerned with the Nori fundamental gerbe rather than the étale one, and hence the theorems regarding the étale fundamental gerbe are not expressed in the right generality. In particular, they always work with inflexible fibered categories, while geometrically connected is the right hypothesis. See also \cite[§ 2,3,4]{tz19}, where part of what is contained in this appendix is done under minor additional hypotheses.

In addition to putting Borne and Vistoli's work for étale fundamental gerbes in the right generality, we give proofs of two technical facts, i.e. the fact that in characteristic $0$ the étale fundamental gerbe behaves well with respect to any field extension (while in \cite{bv15} only algebraic extensions are treated) and the behaviour of the étale fundamental gerbe under finite étale covers. Again, these are not original ideas, but no proof of them was available in the literature.

We want to stress that our effort to state results in maximal generality is not for its own sake: it just happens to work with rather nasty objects that are not even algebraic stacks, like the infinite root stacks of \autoref{open}. Since the theory works for raw fibered categories without any additional hypothesis, we want to give statements in this generality.

If we say that a stack $X$ is finite over a field $k$, we mean that it has groupoid presentation $R\rightrightarrows U$ with $R,U$ finite over $k$, see \cite[§ 4]{bv15}. In particular, a finite stack over $k$ is not necessarily representable.

\subsection{Connected fibered categories}

\begin{definition} \cite[Definition 2.5]{tz19}
	A fibered category $X$ over $k$ is connected if $\H^{0}(X,\O_{X})$ has no nontrivial idempotents.
\end{definition}

\begin{definition}
If $S$ is a scheme and $X$ is a fibered category, we say that a morphism $X\to S$ is \emph{set-surjective} if for every point $s\in S$ there exist a field $\Omega$ and a morphism $\spec \Omega\to X$ with image $s$ in $S$.
\end{definition}

\begin{lemma}
	A fibered category $X/k$ is not connected if and only there exists a set-surjective morphism $X\to\spec k\sqcup\spec k$.
	\begin{proof}
		If $X\to\spec k\sqcup\spec k=\spec k\times k$ is set-surjective, the pullback of $1\times 0$ is a nontrivial idempotent. On the other hand, if $e\in\H^{0}(X,\O_{X})$ is a nontrivial idempotent and $S\to X$ is a morphism, we can define a morphism $S\to\spec k\sqcup\spec k$ by sending $S_{e=0}$ to one point and $S_{e=1}$ to the other one. This defines a morphism $X\to\spec k\sqcup\spec k$. Since $e$ is nontrivial, then for some schemes $S,S'$ with morphisms $S,S'\to X$ we have $S_{e=0}\neq \emptyset$ and $S'_{e=1}\neq\emptyset$, i.e. $X\to\spec k\sqcup\spec k$ is set-surjective.
	\end{proof}
\end{lemma}

Let $X_{1}$, $X_{2}$ be two fibered categories over $k$. It is possible to define the disjoint union $X_{1}\sqcup X_{2}$: is $S$ is a scheme, a morphism $S\to X_{1}\sqcup X_{2}$ is a decomposition of $S=S_{1}\sqcup S_{2}$ with $S_{1},S_{2}$ open and closed together with a pair of morphisms $s_{i}:S_{i}\to X_{i}$.

\begin{definition}
	We define the \emph{clopen} topology on the category of schemes as the Grothendieck topology for which a cover $\{U_{i}\to U\}_{i}$ is a jointly surjective set of morphisms $U_{i}\to U$ which are both closed and open immersions.
\end{definition}

The clopen topology is very coarse, in particular is coarser than the Zariski topology.

\begin{lemma}
	If $X$ is a connected fibered category over $k$ and $X\simeq X_{1}\sqcup X_{2}$, then either $X_{1}$ or $X_{2}$ is empty. If $X$ is a stack in the clopen topology the converse hold, i.e. we can write it as a non-trivial disjoint union if and only if it is disconnected.
	\begin{proof}
		If $X_{1}$ and $X_{2}$ are both non-empty, $1\times 0$ in $\H^{0}(X,\O_{X})=\H^{0}(X_{1},\O_{X_{1}})\times\H^{0}(X_{2},\O_{X_{2}})$ is a nontrivial idempotent.
		
		Let now $e\in\H^{0}(X,\O_{X})$ be a nontrivial idempotent. For every scheme $S$ define
		\[X_{1}(S)=\{s\in X(S)| s^{*}e=1\in\H^{0}(S,\O_{S})\},\]
		\[X_{2}(S)=\{s\in X(S)| s^{*}e=0\in\H^{0}(S,\O_{S})\}.\]
		We have a natural morphism $X\to X_{1}\sqcup X_{2}$ sending a morphism $s:S\to X$ to the pair $s_{1},s_{2}$ where $s_{1}$ is the restriction of $s$ to $S_{e=1}$ and $s_{2}$ is the restriction of $s$ to $S_{e=0}$. Since $S_{e=1}$ and $S_{e=0}$ are open subsets of $S$ such that $S_{e=0}\sqcup S_{e=1}=S$, if $X$ is a stack in the clopen topology we get that $X\to X_{1}\sqcup X_{2}$ is an equivalence.
	\end{proof}
\end{lemma}

\begin{remark}
	If $X$ is an algebraic stack, this is equivalent to asking that the underlying topological space $|X|$ (see \cite[\href{http://stacks.math.columbia.edu/tag/04XE}{Tag 04XE}]{stacks-project}) is connected. On one hand, if $X=X_{1}\sqcup X_{2}$, then $|X|=|X_{1}|\sqcup|X_{2}|$. On the other hand, if $|X|=U_{1}\sqcup U_{2}$ is disconnected, the fact that for every scheme $S$ the natural morphism $|S|\to|X|$ is continuous allows us to define two fibered categories $X_{1}, X_{2}$ such that $|X_{i}|=U_{i}$ and $X=X_{1}\sqcup X_{2}$.
\end{remark}

\subsection{Geometrically connected fibered categories}

If $k'/k$ is a finite extension of fields, the Weil restriction along $k'/k$ is the right adjoint to the functor of base change along $\spec k'\to \spec k$. More concretely, if $X$ is a fibered category over $k$ and $Y$ is a fibered category over $k'$, the Weil restriction $\on{R}_{k'/k}Y$ is a fibered category over $k$ with an equivalence of categories
\[\hom_{k}(X,\on{R}_{k'/k}Y)\simeq\hom_{k'}(X_{k'},Y)\]
functorial in $X$ and $Y$. We can construct $\on{R}_{k'/k}Y$ as the fibered product $\aff/k\times_{\aff/k'}Y$. When $Y$ is represented by a quasi-projective scheme, $\on{R}_{k'/k}Y$ is represented a scheme, too. If $Y$ is represented by a finite stack and $k'/k$ is separable, then $\on{R}_{k'/k}Y$ is represented by a finite stack too, see \cite[Lemma 6.2]{bv15}.

\begin{lemma}\label{weil}
	Let $k'/k$ be a finite, separable extension, and $Y$ a finite étale stack over $k'$. Then $\on{R}_{k'/k}Y$ is a finite étale stack over $k$, too.
	\begin{proof}
		In the proof of \cite[Lemma 6.2]{bv15}, from a finite groupoid presentation $R\rightrightarrows U$ of $Y$ they construct a finite groupoid presentation $R'\rightrightarrows U'$ of $\on{R}_{k'/k}Y$. Following their construction, it is immediate to check that if $R\rightrightarrows U$ is étale, $R'\rightrightarrows U'$ is étale too.
	\end{proof}
\end{lemma}

Recall that a fibered category is \emph{concentrated} if there exists an affine scheme $U$ and a representable, quasi separated, quasi compact and faithfully flat morphism $U\to X$.

If $X$ is concentrated and $u:U\to X$ is as above, set $R=U\times_{X}U$, we obtain an fpqc groupoid $(r_{1},r_{2}):R\rightrightarrows U$ in algebraic spaces. From standard arguments in descent theory we get an exact sequence
\[0\to\H^{0}(X,\O_{X})\xar{u^{*}}\H^{0}(U,\O_{U})\xar{r_{1}^{*}-r_{2}^{*}}\H^{0}(R,\O_{R})\]
and hence it follows easily that for any field extension $k'/k$,
\[\H^{0}(X_{k'},\O_{X_{k'}})=\H^{0}(X,\O_{X})\otimes_{k}k'.\]

\begin{lemma}\label{geomconn}
	Let $X$ be a category fibered over $k$, and $k_{s}/k$ a separable closure. Consider the following:
	\begin{enumerate}[(i)]
		\item $X_{k'}/k'$ is connected for every extension $k'/k$,
		\item $X_{k_{s}}/k_{s}$ is connected,
		\item $X_{k'}/k'$ is connected for every finite, separable extension $k'/k$,
		\item $k$ is the only étale subalgebra of $\H^{0}(X,\O_{X})$,
		\item $\spec\H^{0}(X,\O_{X})$ is geometrically connected.
	\end{enumerate}
	In general, we have implications $(i)\Leftrightarrow (ii)\Rightarrow (iii)\Leftrightarrow (iv)\Leftrightarrow(v)$. If $X$ is an algebraic space or it is concentrated, then $(iii)\Rightarrow (ii)$ holds, too.

	\begin{proof}
		\begin{description}
			\item[$(i)\Rightarrow (ii)$] Obvious.
			
			\item[$(ii)\Rightarrow (i)$] Suppose that $X_{k'}\to\spec k'\sqcup\spec k'$ is a set-surjective morphism. Up to enlarging $k'$, we may suppose that $k_{s}\s k'$. Let $S$ be a scheme over $k_{s}$, and $S\to X_{k_{s}}$ a morphism. By \cite[\href{http://stacks.math.columbia.edu/tag/0363}{Tag 0363}]{stacks-project} and \cite[\href{http://stacks.math.columbia.edu/tag/0383}{Tag 0383}]{stacks-project}, $S_{k'}\to S$ is open and induces a bijection of connected components.
			
			In particular, we can write $S=S_{1}\sqcup S_{2}$ such that $S_{i,k'}\to S_{k'}\to\spec k'\sqcup\spec k'$ maps to the $i$-th point, for $i=1,2$. This allows us to define a morphism $S\to\spec k_{s}\sqcup\spec k_{s}$ whose base change is $S_{k'}\to\spec k'\sqcup\spec k'$, and thus a morphism $X_{k_{s}}\to\spec k_{s}\sqcup\spec k_{s}$ whose base change is $X_{k'}\to\spec k'\sqcup\spec k'$. The morphism $X_{k_{s}}\to\spec k_{s}\sqcup\spec k_{s}$ is clearly set-surjective, and this is absurd.
			
			\item[$(ii)\Rightarrow (iii)$] If $X_{k'}\to\spec k'\sqcup\spec k'$ is set-surjective, then $X_{k_{s}}\to\spec k_{s}\sqcup\spec k_{s}$ is set-surjective too.
			
			\item[$(iii)\Rightarrow (iv)$] Suppose that $A\s\H^{0}(X,\O_{X})$ is a nontrivial finite étale subalgebra of degree $d>1$, there exists a scheme $S$ with a morphism $S\to X$ such that the composition $S\to X\to \spec A$ is dominant. Now choose $k'/k$ a finite separable extension which splits $A$. The base change
			\[X_{k'}\to\spec A_{k'}=\spec k'^{d}\]
			is set-surjective because $S_{k'}\to X_{k'}\to\spec k'^{d}$ is set-surjective. But this is absurd, since $d>1$ and $X_{k'}$ is connected.
			
			\item[$(iv)\Rightarrow (iii)$] Suppose that $k'/k$ is a finite separable extension and that we have a set-surjective morphism $X_{k'}\to\spec k'\sqcup\spec k'$, this induces a morphism $X\to\on{R}_{k'/k}(\spec k'\sqcup\spec k')$. Since $\on{R}_{k'/k}(\spec k'\sqcup\spec k')$ is a finite étale scheme, by hypothesis we have a factorization
			\[X\to\spec k\to\on{R}_{k'/k}(\spec k'\sqcup\spec k').\]
			But this gives a factorization
			\[X_{k'}\to\spec k'\to\spec k'\sqcup\spec k'\]
			which is absurd.
			
			\item[$(iv)\Leftrightarrow (v)$] This is well known.
		\end{description}
		
		For the implication $(iii)\Rightarrow (ii)$, if $X$ is concentrated we have
		\[\H^{0}(X_{k'},\O_{X_{k'}})=\H^{0}(X,\O_{X})\otimes_{k}k'\]
		for every extension $k'/k$, hence we can reduce to affine schemes for which the result is well known. If $X$ is an algebraic space, this is \cite[\href{https://stacks.math.columbia.edu/tag/0A17}{Lemma 0A17}]{stacks-project}.
	\end{proof}
\end{lemma}

\begin{definition}\label{geomconndef}
	Let $X$ be a fibered category. We say that $X$ is \emph{geometrically connected} if the equivalent conditions $(iii)$, $(iv)$ and $(v)$ of \autoref{geomconn} hold for $X$.
\end{definition}

\subsection{Existence and base change}

\begin{lemma}\label{fetstack}
	A stack over $k$ is finite étale if and only if it is étale and of finite type over $k$.
	\begin{proof}
		A finite stack in the sense of \cite[§ 4]{bv15} is clearly proper. On the other hand, let $X$ be an étale stack of finite type over $k$. Étale morphisms are by definition DM, thus $X$ is a DM stack. Let $U\to X$ be an étale cover of finite type with $U$ a scheme and $R=U\times_{R}U$. By composition, we get that $U,R$ are étale of finite type over $k$. It follows that $X$ is finite étale.
	\end{proof}
\end{lemma}

\begin{definition}
	An fpqc stack $\Gamma$ over a field $k$ is \emph{profinite étale} if it is the limit of a projective system of finite, étale stacks over $k$, in the sense of \cite[Definition 3.5]{bv15}.
\end{definition}

\begin{remark}
	In \cite[Definition 3.5]{bv15} they define the limit of a projective system $(\Gamma_{i})_{i}$ of affine fpqc gerbes as a category fibered in groupoids which turns out to be an fpqc stack. Actually, it is straightforward to check that the definition works without any modification for a projective system $(\Gamma_{i})$ of categories fibered in groupoids, and if $\Gamma_{i}$ is an fpqc stack for every $i$ then also the limit is an fpqc stack. Moreover, if $\Gamma_{i}$ is an affine fpqc gerbe for every $i$ and the limit is not empty, then the limit is an fpqc gerbe too, see \cite[Proposition 3.7]{bv15}.
\end{remark}

\begin{definition}
	Let $X$ be a fibered category over $k$, and $\Pi$ a profinite étale gerbe with a morphism $X\to\Pi$. Then $X\to\Pi$ is an \emph{étale fundamental gerbe} if, for every finite, étale stack $\Phi$, the functor
	\[\hom(\Pi,\Phi)\to\hom(X,\Phi)\]
	is an equivalence of categories.
\end{definition}

\begin{lemma}\label{proetale}
	Let $X$ be a fibered category with an étale fundamental gerbe $X\to\Pi$, and $\Phi$ a profinite étale stack. Then
	\[\hom(\Pi,\Phi)\to\hom(X,\Phi)\]
	is an equivalence of categories.
	In particular, the étale fundamental gerbe is unique up to a canonical equivalence.
	\begin{proof}
		This is a straightforward application of the definition of the étale fundamental gerbe and of profinite étale stacks.
	\end{proof}
\end{lemma}

The following simple lemma is rather enlightening in the sense that it draws the line between the étale setting and the Nori setting: its failure for finite stacks is what makes Nori's fundamental gerbe subtler than the étale one. 

\begin{lemma}\label{etalegerbe}
	Let $\Phi$ be a finite étale stack. Then the natural morphism 
	\[\Phi\to\spec\H^{0}(\Phi,\O_{\Phi})\]
	is a gerbe.
	\begin{proof}
		We give an elementary proof. See also \cite[Proposition 3.2]{tz19} for a more technical proof for finite, reduced stacks.
		
		If $k^{s}/k$ is the separable closure, it is easy to check that $\Phi\to\spec\H^{0}(\Phi,\O_{\Phi})$ is a gerbe if and only if $\Phi_{k^{s}}\to\spec\H^{0}(\Phi_{k^{s}},\O_{\Phi_{k^{s}}})$ is a gerbe. Hence, we may suppose that $k$ is separably closed.
		
		Choose now a finite étale groupoid $R\rightrightarrows U$ giving a presentation of $\Phi$. Since $k$ is separably closed and $R$, $U$ are finite étale, they are simply finite disjoint unions of points. Hence we can write
		\[\Phi=\sqcup_{i}BG_{i}\]
		where $G_{i}$ are finite discrete groups. Now it is obvious that
		\[\Phi=\sqcup_{i}BG_{i}\to\sqcup_{i}\spec k\]
		is a gerbe.
	\end{proof}
\end{lemma}

\begin{corollary}\label{etaleinflexible}
	Let $X$ be a fibered category. Then $X$ is geometrically connected if and only if every morphism $X\to\Gamma$ where $\Gamma$ is a finite étale stack has a factorization
	\[X\to\Gamma'\to\Gamma\]
	where $\Gamma'$ is a finite étale gerbe.
	\begin{proof}
		Suppose that $X$ is geometrically connected. Consider the composition
		\[X\to\Gamma\to\spec\H^{0}(\Gamma,\O_{\Gamma}).\]
		Since $X$ is geometrically connected and $\H^{0}(\Gamma,\O_{\Gamma})$ is finite étale, we have a factorization
		\[X\to\spec k\to\spec\H^{0}(\Gamma,\O_{\Gamma}).\]
		Set $\Gamma'=\spec k\times_{\spec\H^{0}(\Gamma,\O_{\Gamma})}\Gamma$, we have a factorization
		\[X\to\Gamma'\to\Gamma\]
		and $\Gamma'$ is a gerbe over $\spec k$ thanks to \autoref{etalegerbe}.
		
		On the other hand, if $A\s\H^{0}(X,\O_{X})$ is a nontrivial étale subalgebra, the natural morphism $X\to\spec A$ cannot factorize through any finite gerbe.
	\end{proof}
\end{corollary}

The following three results are straightforward modifications of results of Borne and Vistoli.

\begin{theorem}[{\cite[Theorem 5.7]{bv15}}]\label{existence}
	Let $X$ be a fibered category over $k$. Then $X$ has an étale fundamental gerbe if and only if it is geometrically connected.
	\begin{proof}
		Thanks to \autoref{etaleinflexible}, we can replace inflexible fibered categories with geometrically connected ones. See also \cite[Proposition 4.3]{tz19} for a proof under some minor additional hypotheses.
	\end{proof}
\end{theorem}

\begin{proposition}[{\cite[Proposition 6.1]{bv15}}]\label{basechangefin}
	Let $k'/k$ be an algebraic and separable extension, $X$ a geometrically connected fibered category over $k$. Suppose that either
	\begin{enumerate}[(a)]
		\item $k'$ is finite over $k$, or
		\item $X$ is concentrated.
	\end{enumerate}
	Then $X_{k'}$ is geometrically connected over $k'$ and $\Pi_{X_{k'}/k'}=\spec k'\times\Pi_{X/k}$.
	\begin{proof}
		Replace \cite[Lemma 6.2]{bv15} with \autoref{weil}.
	\end{proof}
\end{proposition}

Given two extensions $G,H$ of a group $\Gamma$, we have defined $\on{Hom-ext}_{\Gamma}(G,H)$ as the set of homomorphisms $G\to H$ of $\Gamma$-extensions modulo the action of $\ker(H\to\Gamma)$ by conjugacy. It is more natural to consider the category $\on{\bf Hom-ext}_{\Gamma}(G,H)$ whose objects are said homomorphisms and whose arrows are given by conjugacy with elements of $\ker(H\to\Gamma)$, then $\on{Hom-ext}_{\Gamma}(G,H)$ is the set of isomorphism classes.

\begin{proposition}[{\cite[Proposition 9.3]{bv15}}]\label{fundsect}
	Let $X$ be a quasi-compact, quasi-separated and geometrically connected algebraic stack over $k$ with a geometric point $\bar{x}:\spec\Omega\to X$, and $T$ any geometrically connected scheme with a geometric point $\bar{t}:\spec\Omega\to X$. There is a (non-canonical) equivalence of categories
	\[\Pi_{X/k}(T)\to\on{\bf Hom-ext}_{\Gamma_{k}}(\pi_{1}(T,\bar{t}),\pi_{1}(X,\bar{x}))\]
	that composed with the canonical functor $\hom_{k}(T,X)\to\Pi_{X/k}(T)$ gives the natural map 
	\[\hom_{k}(T,X)\to\on{\bf Hom-ext}_{\Gamma_{k}}(\pi_{1}(T,\bar{t}),\pi_{1}(X,\bar{x})).\]
	\begin{proof}
		Replacing $\spec k$ with $T$ simply doesn't affect the proof.
	\end{proof}
\end{proposition}

Suppose now that we are in characteristic $0$. Following Borne and Vistoli, we have shown that the étale fundamental gerbe behaves well under algebraic field extensions: we want to show that, actually, it behaves well with respect to any field extension. The idea is to rephrase the theorem in terms of étale fundamental groups, and then use the fact that the étale fundamental group is invariant along extensions of algebraically closed fields, see \cite[Exposé XIII, Proposition 4.6]{sga1}.

Given a field $k$, let $\affgrp$ be the category of affine group schemes over $k$, and let $\fgrp$, $\on{PFGrp}$ be the categories of finite and profinite classical groups.

\begin{lemma}
	The functor $\fgrp\to\affgrp$ which sends a finite group to the associated constant group scheme extends to a fully faithful functor $\on{PFGrp}\to\affgrp$ which preserves limits.
	
	If $k$ is separably closed, the essential image of the functor $\on{PFGrp}\to\affgrp$ consists of the full subcategory of profinite étale group schemes.
	\begin{proof}
		Given a profinite group $G$, define $I_{G}$ as the $k$-algebra of \emph{continuous} homomorphisms $G\to k$ where $k$ is endowed with the discrete topology. The $k$-algebra $I_{G}$ can be endowed with an Hopf algebra structure analogously to the finite case. Define the associated profinite constant group scheme
		\[\underline{G}=\spec I_{G}.\]
		If we write $G=\projlim_{i}G_{i}$ as a limit of finite groups, the natural map $\injlim_{i}I_{G_{i}}\to I_{G}$ is easily checked to be an isomorphism of Hopf algebras. The rest of the proof follows directly.
	\end{proof}
\end{lemma}

\begin{lemma}\label{pfgrp}
	If $G,H$ are profinite étale groups schemes and $k'/k$ is an extension of separably closed fields, then the natural functor
	\[\hom_{k}(B_{k}G,B_{k}H)\to\hom_{k'}(B_{k'}G,B_{k'}H)\]
	is an equivalence.
	\begin{proof}
		Both categories have the same description in purely group theoretic terms: let us explain this. Since $k$ is separably closed, thanks to \autoref{pfgrp} there exist profinite groups $G_{c},H_{c}$ whose associated group schemes over $k$ are $G,H$. Clearly, $G_{k'},H_{k'}$ are associated to $G_{c},H_{c}$ over $k'$.
		
		We define now the category $\hom_{\rm cont}(BG_{c},BH_{c})$ of "continuous" functors $BG_{c}\to BH_{c}$: its objects are just continuous homomorphisms $G_{c}\to H_{c}$, and every $h\in H_{c}$ defines an arrow $\phi\to h^{-1}\phi h$ for every continuous homomorphism $\phi:G_{c}\to H_{c}$.
		
		We have a natural morphism $\hom_{\rm cont}(BG_{c},BH_{c})\to\hom(B_{k}G,B_{k}H)$: since $k$ is separably closed, thanks to \autoref{pfgrp} it is immediate to check that it is an equivalence of categories. The same is true over $k'$, thus the thesis follows by 2-commutativity of the obvious diagram.
	\end{proof}
\end{lemma}

For the following \autoref{fund}, I would like to thank Marc Hoyois who suggested the use of noetherian approximation in order to reach full generality, see \href{https://mathoverflow.net/questions/294847/}{MathOverflow 294847}.

\begin{lemma}\label{fund}
	Let $k'/k$ be an extension of algebraically closed fields of characteristic $0$. Consider $X$ a concentrated fibered category over $k$, and $\Phi$ a finite étale stack over $k$. Then the natural functor
	\[\hom_{k}(X,\Phi)\to\hom_{k'}(X_{k'},\Phi_{k'})\]
	is an equivalence of categories.
	\begin{proof}
		Let us prove this firstly under the additional hypothesis that $X$ is a scheme of finite type over $k$. Under this hypothesis, connected components are open, hence we may suppose that $X$ is connected and $\Phi$ is of the form $BG$ for some finite group $G$. Fix any point $x\in X(k)$. Thanks to \cite[Exposé XIII, Proposition 4.6]{sga1}, $\pi_{1}(X,x)=\pi_{1}(X_{k'},x_{k'})$.
		
		We have thus
		\[\hom_{k}(X,B_{k}G)=\hom_{k}(B_{k}\pi_{1}(X,x),B_{k}G)=\]
		\[=\hom_{k'}(B_{k'}\pi_{1}(X_{k'},s_{k'}),B_{k'}G)=\hom_{k'}(X_{k'},B_{k'}G).\]
		
		Let us now generalize to $X$ quasi compact, quasi separated scheme.	By noetherian approximation \cite[Theorem C.9]{tt90}, we can write $X$ as an inverse limit $\projlim_{i}X_{i}$, with $X_{i}$ of finite type over $k$. Since $\Phi$ is finite,
		\[\hom_{k}(X,\Phi)=\injlim_{i}\hom_{k}(X_{i},\Phi)=\]
		\[=\injlim_{i}\hom_{k'}(X_{i,k'},\Phi_{k'})=\hom_{k'}(X_{k'},\Phi_{k'}).\]
		
		Finally, if $X$ is a concentrated fibered category, let $U$ be a quasi compact and quasi separated scheme with a representable, quasi separated, quasi compact and faithfully flat morphism $U\to X$. Set $R=U\times_{X}U$, $R$ is again quasi compact and quasi separated. Let $\hom(R\rra U,\Phi)$ be the category of morphism $U\to\Phi$ satisfying the usual cocycle condition on $R$. Descent theory tells us that $\hom(R\rra U,\Phi)$ is naturally equivalent to $\hom(X,\Phi)$, even if $X$ is not a stack and hence $X\neq [U/R]$. Since $U$ and $R$ are quasi-compact and quasi separated, by the preceding case we conclude that
		\[\hom_{k'}(X_{k'},\Phi_{k'})=\hom_{k'}(R_{k'}\rra U_{k'},\Phi_{k'})=\]
		\[=\hom_{k}(R\rra U,\Phi)=\hom_{k}(X,\Phi).\]
	\end{proof}
\end{lemma}

\begin{proposition}\label{basechange}
	Let $k$ be a field of characteristic $0$. If $X$ is a geometrically connected, concentrated fibered category over $k$, then the natural map $\Pi_{X_{k'}/k'}\to\Pi_{X/k}\times_{k}k'$ is an isomorphism for every field extension $k'/k$.
	\begin{proof}
		Thanks to \autoref{basechangefin}, it is immediate to reduce to the case in which $k$ and $k'$ are both algebraically closed. We have to show that $\Pi_{X/k}\times_{k}k'$ has the universal property of the étale fundamental gerbe of $X_{k'}$.
		
		Since $k'$ is algebraically closed, every finite étale stack over $k'$ has the form $\sqcup_{i}B_{k'}G_{i}$ for some finite number of finite groups $G_{i}$. In particular, every finite étale stack over $k'$ is isomorphic to $\Phi_{k'}$ for some finite étale stack $\Phi$ over $k$, hence it is enough to show that $\Pi_{X/k}\times_{k}k'$ has the universal property with respect to stacks of the form $\Phi_{k'}$ with $\Phi$ finite étale over $k$.
		
		Now observe that $\Pi_{X/k}$, being a gerbe over $\spec k$, is concentrated: in fact, any morphism $\spec L\to\Pi_{X/k}$ with $L$ a field is representable, quasi compact, quasi separated and faithfully flat. Hence both $X$ and $\Pi_{X/k}$ are concentrated and thanks to \autoref{fund} we have	
		\[\hom_{k'}(X_{k'},\Phi_{k'})=\hom_{k}(X,\Phi)=\]
		\[=\hom_{k}(\Pi_{X/k},\Phi)=\hom_{k'}(\Pi_{X/k}\times_{k}k',\Phi_{k'}).\]
	\end{proof}
\end{proposition}

\subsection{Étale coverings of fibered categories}

\begin{lemma}\label{degree}
	Let $f:Y\to X$ be a representable, finite étale morphism of fibered categories. If $X$ is connected then there exists an integer $d$ such that for every scheme $S$ and every morphism $s:S\to X$ the étale covering $S\times_{X}Y\to S$ has constant degree $d$.
	\begin{proof}
		If $S$ is a scheme, $s\in X(S)$ an object and $d\geq 0$ an integer, the locus $S_{=d}$ of points $p$ of $S$ such that $Y\times_{X}S\to S$ has degree $d$ over $p$ is an open and closed sub-scheme of $S$, set $S_{\neq d}=S \setminus S_{=d}$. This allows to define a morphism $X\to\spec k\sqcup\spec k$ sending $S_{=d}$ to the first point and $S_{\neq d}$ to the second point, and if there exist morphisms $S,S'\to X$ such that $S_{=d}$ and $S'_{\neq d}$ are both nonempty then $X$ is not connected, and this is absurd.
		
		There exists some $d_{0}$ and a morphism $S\to X$ such that $S_{=d_{0}}\neq\emptyset$, hence for every morphism $S'\to X$ we have $S'_{=d_{0}}=S'$, i.e. $Y\to X$ has constant degree $d_{0}$.
	\end{proof}
\end{lemma}

In the following, we need to use quotients by group actions of stacks, see \cite[Theorem 4.1]{rom05} for the general existence result.

\begin{proposition}\label{cartetale}
	Let $Y\to X$ be a representable, finite étale morphism of geometrically connected fibered categories. The following natural $2$-commutative diagram is $2$-cartesian.
	\[\begin{tikzcd}
		Y\rar\dar	&	\Pi_{Y/k}\dar	\\
		X\rar		&	\Pi_{X/k}
	\end{tikzcd}\]
\begin{proof}
	Thanks to \autoref{degree}, $Y\to X$ is a finite cover of fixed degree $d$. Let $d\times X$ be the disjoint union of $d$ copies of $X$, we have a finite cover $d\times X\to X$ of degree $d$. The group $S_{d}$ acts on the fibered category $Z=\uisom_{X}(d\times X,Y)$ by automorphisms of $d\times X$ making it into an $S_{d}$-torsor over $X$. If $S$ is a scheme with a morphism $S\to Z$, we have a trivialization $d\times S\simeq Y\times_{X}S$. The first copy of $d\times S$ gives us a morphism $S\to Y$, and thus by Yoneda's lemma we have a $S_{d-1}$-invariant morphism $Z\to Y$ which is actually a $S_{d-1}$-torsor.
	
	All of this can be packed by saying that we have a morphism $X\to BS_{d}$ with identifications $Z=X\times_{BS_{d}}\spec k$ and $Y=X\times_{BS_{d}}BS_{d-1}$. Moreover, define $\Pi=\Pi_{X/k}\times_{BS_{d}}BS_{d-1}$ and $\Lambda=\spec k\times_{BS_{d}}\Pi_{X/k}$. We have a $2$-cartesian diagram
	\[\begin{tikzcd}
		Z\rar\dar	&	\Lambda\rar\dar		&	\spec k\dar		\\
		Y\rar\dar	&	\Pi\rar\dar			&	BS_{d-1}\dar	\\
		X\rar		&	\Pi_{X/k}\rar		&	BS_{d}
	\end{tikzcd}\]
	
	Since $\Pi$ is profinite étale, if we show that it satisfies the universal property of $\Pi_{Y/k}$ then we have that $\Pi=\Pi_{Y/k}$ thanks to \autoref{proetale}, and hence the thesis.	
	
	Consider now a finite étale stack $\Phi$: we want to show that
	\[\hom_{k}\left(\Pi,\Phi\right)\to\hom_{k}\left(Y,\Phi\right)\]
	is an equivalence of categories. Let $\rho:Z\times S_{d}\to Z$ be the action. If $Y\to\Phi$ is a morphism, consider the composition
	\[\rho_{\Phi}:Z\times S_{d}\xar{\rho}Z\to Y\to\Phi.\]
	For every $g\in S_{d}$, this defines a morphism $\rho_{\Phi}(\cdot,g):Z\to\Phi$. If $h\in S_{d-1}\s S_{d}$, since $Z\to Y$ is $S_{d-1}$ invariant we get that $\rho_{\Phi}(\cdot,g)=\rho_{\Phi}(\cdot,gh):Z\to\Phi$, hence $\rho_{\Phi}(\cdot,[g])$ is well defined for $[g]\in S_{d}/S_{d-1}$. This gives us an $S_{d}$-equivariant morphism
	\[Z\to\Phi^{S_{d}/S_{d-1}}\]
	where $S_{d}$ acts on $\Phi^{S_{d}/S_{d-1}}$ via left multiplication on $S_{d}/S_{d-1}$.

	On the other hand, if we have an $S_{d}$-equivariant morphism $Z\to\Phi^{S_{d}/S_{d-1}}$, it is $S_{d-1}$-invariant since $S_{d-1}$ acts trivially on $S_{d}/S_{d-1}$. Hence we have an induced morphism
	\[Y\to\Phi^{S_{d}/S_{d-1}},\]
	which, composed with the projection $\Phi^{S_{d}/S_{d-1}}\to\Phi$ on the identity component, gives a morphism $Y\to\Phi$. It is easy to check that these constructions are inverses and give an equivalence of categories
	\[\hom(Y,\Phi)\xar{\sim}\hom^{S_{d}}(Z,\Phi^{S_{d}/S_{d-1}}).\]
	Since $Z\to X$ is an $S_{d}$-torsor, we also have an equivalence
	\[\hom^{S_{d}}(Z,\Phi^{S_{d}/S_{d-1}})\xar{\sim}\hom_{BS_{d}}(X,[\Phi^{S_{d}/S_{d-1}}/S_{d}])\]
	and their composition
	\[\hom(Y,\Phi)\xar{\sim}\hom_{BS_{d}}(X,[\Phi^{S_{d}/S_{d-1}}/S_{d}]).\]
	We can repeat the same argument with $\Pi_{X/k}$, $\Pi$ and $\Lambda$ instead of $X$, $Y$ and $Z$, finding an equivalence
	\[\hom(\Pi,\Phi)\xar{\sim}\hom_{BS_{d}}(\Pi_{X/k},[\Phi^{S_{d}/S_{d-1}}/S_{d}]).\]
	Thanks to \autoref{fetstack}, $[\Phi^{S_{d}/S_{d-1}}/S_{d}]$ is a finite étale stack and thus there is another equivalence
	\[\hom_{BS_{d}}(X,[\Phi^{S_{d}/S_{d-1}}/S_{d}])\xar{\sim}\hom_{BS_{d}}(\Pi_{X/k},[\Phi^{S_{d}/S_{d-1}}/S_{d}]).\]
	Composing these three, we obtain the desired equivalence
	\[\hom(Y,\Phi)\xar{\sim}\hom(\Pi,\Phi).\]
\end{proof}
\end{proposition}

\section*{Acknowledgements}

This article is part of my PhD thesis. I would like to thank my PhD advisor Angelo Vistoli for many useful discussions and for teaching me how to use stacks (and almost everything else). I would also like to thank Tamás Szamuely and Hélène Esnault for many useful remarks, and Marc Hoyois for pointing me out how to use noetherian approximation in order to reach full generality in \autoref{fund}, see \href{https://mathoverflow.net/questions/294847/}{MathOverflow 294847}. Finally, I would like to thank an anonymous referee for carefully reading an earlier draft of the paper and helping me to greatly clarify the exposition.

\bibliographystyle{amsalpha}
\bibliography{anab}

\end{document}